\pdfoutput=1
\RequirePackage{ifpdf}
\ifpdf % We are running pdfTeX in pdf mode
\documentclass[pdftex]{sigma}
\else
\documentclass{sigma}
\fi

\numberwithin{equation}{section}

{ \theoremstyle{definition}
\newtheorem{dfn}{Definition}[section]

\newtheorem{rmrk}[dfn]{Remark}
}
\newtheorem{thm}[dfn]{Theorem}
\newtheorem{lmma}[dfn]{Lemma}
\newtheorem{ppsn}[dfn]{Proposition}
\newtheorem{crlre}[dfn]{Corollary}
\newtheorem*{theorem-nonumber}{Theorem}

\begin{document}

\newcommand{\arXivNumber}{1307.4850}

\allowdisplaybreaks

\renewcommand{\thefootnote}{$\star$}

\renewcommand{\PaperNumber}{076}

\FirstPageHeading

\ShortArticleName{Quantum Isometry Groups of Noncommutative Manifolds}

\ArticleName{Quantum Isometry Groups of Noncommutative\\
Manifolds Obtained by Deformation Using Dual\\
Unitary 2-Cocycles\footnote{This paper is a~contribution to the Special Issue on Noncommutative Geometry and Quantum
Groups in honor of Marc A.~Rief\/fel.
The full collection is available at
\href{http://www.emis.de/journals/SIGMA/Rieffel.html}{http://www.emis.de/journals/SIGMA/Rieffel.html}}}

\Author{Debashish GOSWAMI and Soumalya JOARDAR}

\AuthorNameForHeading{D.~Goswami and S.~Joardar}

\Address{Indian Statistical Institute, 203, B.T.~Road, Kolkata 700108, India}
\Email{\href{mailto:goswamid@isical.ac.in}{goswamid@isical.ac.in},
\href{mailto:soumalya.j@gmail.com}{soumalya.j@gmail.com}}

\ArticleDates{Received January 29, 2014, in f\/inal form July 11, 2014; Published online July 17, 2014}

\Abstract{It is proved that the (volume and orientation-preserving) quantum isometry group of a~spectral triple obtained
by deformation by some dual unitary 2-cocycle is isomorphic with a~similar twist-deformation of the quantum isometry
group of the original (undeformed) spectral triple.
This result generalizes similar work by Bhowmick and Goswami for Rief\/fel-deformed spectral triples in~[\textit{Comm. Math. Phys.} \textbf{285} (2009), 421--444].}


 \Classification{58B34; 46L65; 81R50}

\rightline{\it Dedicated to Professor Marc A.~Rieffel}

\renewcommand{\thefootnote}{\arabic{footnote}}
\setcounter{footnote}{0}

\section{Introduction}

It is a~very important and interesting problem in the theory of quantum groups and noncommutative geometry to study
`quantum symmetries' of various classical and quantum structures.
Indeed, symmetries of physical systems (classical or quantum) were conventionally modeled by group actions, and after
the advent of quantum groups, group symmetries were naturally ge\-ne\-ra\-lized to symmetries given by quantum group actions.
In this context, it is natural to think of quantum automorphism or the full quantum symmetry groups of various
mathematical and physical structures.
The underlying basic principle of def\/ining a~quantum automorphism group of a~given mathematical structure consists of
two steps: f\/irst, to identify (if possible) the group of automorphisms of the structure as a~universal object in
a~suitable category, and then, try to look for the universal object in a~similar but bigger category by replacing groups
by quantum groups of appropriate type.
Initiated by S.~Wang who def\/ined and studied quantum permutation groups of f\/inite sets and quantum automorphism groups
of f\/inite-dimensional algebras, such questions were taken up by a~number of mathematicians including Banica, Bichon
(see, e.g.,~\cite{ban_2,graph,wang}), and more recently in the framework of Connes' noncommutative
geometry~\cite{con} by Goswami, Bhowmick, Skalski, Soltan, Banica and others who have extensively studied the quantum
group of isometries (or quantum isometry group) def\/ined in~\cite{Goswami} (see
also~\cite{qdisc,qorient,qiso_trans} etc.).

Most of the examples of noncommutative manifolds are obtained by deforming classical spectral triples.
It was shown in~\cite{qiso_comp} that the quantum isometry group of a~Rief\/fel-deformed noncommutative manifold can be
obtained by a~similar deformation (Rief\/fel--Wang, see~\cite{rieffel_wang}) of the quantum isometry group of the original
(undeformed) noncommutative manifold.

In the present paper our goal is to generalize Bhowmick--Goswami's results about Rief\/fel-deformation to any cocycle
twisted spectral triple.
Combining this with the fact (proved in~\cite{rigidity}) that the quantum isometry group of a~classical compact,
connected, Riemannian manifold is the same as the classical isometry group of such manifold (i.e.~there is no genuine
quantum isometry for such manifold), we have been able to compute the quantum isometry group of non commutative
manifolds obtained from classical manifolds using unitary 2-cocycles.

Let $({\cal A}^{\infty},{\cal H},{\cal D})$ be a~spectral triple of compact type,~$R$ a~positive, invertible
operator com\-mu\-ting with ${\cal D}$.
Given a~dual unitary 2-cocycle~$\sigma$ on a~quantum subgroup (say ${\cal Q}$) of the quantum isometry group
${\rm QISO}^{+}_{R}({\cal A}^{\infty},{\cal H},{\cal D})$ (in the sense of~\cite{qorient}), we construct a~spectral triple on
the twisted algebra $({\cal A}^{\infty})^{\sigma}$ on the same Hilbert space ${\cal H}$ with the same Dirac operator
${\cal D}$.
Moreover, we construct a~canonical `twisted' operator $R^{\sigma}$ commuting with ${\cal D}$ using the Haar state of
${\cal Q}$ and the cocycle~$\sigma$.
The main result of this paper (Theorem~\ref{twist_main}) can now be stated:

\begin{theorem-nonumber}
The quantum isometry group ${\rm QISO}^{+}_{R^{\sigma}}(({\cal A}^{\infty})^{\sigma},{\cal H},{\cal D})$
is isomorphic with
${\rm QISO}^{+}_{R}({\cal A}^{\infty},$ ${\cal H},{\cal D})^{\sigma}$, i.e.~the cocycle twist of the quantum isometry group of the
original spectral triple.
\end{theorem-nonumber}

In this paper we have also been able to relax the assumption of existence of a~dense $\ast$-algebra for which
the action is algebraic in case of Rief\/fel deformed manifolds.

Let us mention that after we wrote and uploaded to the internet a~preliminary version of this article, Neshveyev and
Tuset posted the article~\cite{nesh2} in which they formulated a~general theory of deformation of $C^*$-algebras by dual
unitary 2-cocycles on locally compact quantum groups.
Indeed, our set-up for von Neumann algebraic deformation is closely related to the restriction of their framework to the
simpler case of compact quantum groups and is nothing but an alternative description of the von Neumann algebra
generated by the $C^*$-algebraic deformation obtained in~\cite{nesh2}.
However, there are two important subtle dif\/ferences between our set-up and that of~\cite{nesh2}: f\/irst, we work with
universal compact quantum groups rather than those in the reduced form as in~\cite{nesh2}, and second, we do not need
$C^*$-algebraic action of the quantum group.
Moreover, the main emphasis of our article is on deformations of spectral triples and their quantum isometry groups, not
just the study of deformed algebras.
We plan to re-cast some of our results in the framework of~\cite{nesh2} and also study deformed spectral triples more
deeply in a~forthcoming paper.

Let us conclude this section with a~brief note on the notations to be used in the paper.
For a~subset~$A$ of a~vector space~$V$, $\operatorname{Sp} A$ denotes the linear span of~$A$.
We use $\oplus$, $\otimes_{\rm alg}$, and $\overline{\otimes}$ to denote algebraic direct sum, algebraic tensor product and
von Neumann algebraic tensor product respectively (see~\cite{Takesaki} for the details on tensor product).
For a~Hilbert space ${\cal H}$, ${\cal K}({\cal H})$ and ${\cal B}({\cal H})$ will denote the $C^{\ast}$-algebra of
compact operators and the $C^{\ast}$-algebra of bounded operators on ${\cal H}$ respectively.
$\otimes$ stands for the exterior tensor product for two Hilbert bimodules (in particular a~Hilbert space ${\cal H}$ and
a~$C^{\ast}$-algebra ${\cal Q}$) and tensor product of Hilbert spaces.
Spatial tensor products of $C^{\ast}$-algebras will be denoted by $\hat{\otimes}$.
We use the standard Kronecker $\delta_{ij}$ symbol to denote the function which is~$1$ if~$i=j$ and~$0$ if~$i\neq j$.
The left convolution of an element in a~CQG with a~linear functional~$f$ will be denoted by $f\ast a$ ($:=(f\otimes
{\rm id})\Delta(a))$ and the right convolution will be denoted by $a\ast f$ ($:=({\rm id}\otimes f)\Delta(a))$ for $a\in{\cal
Q}$, $\Delta$~is the coproduct of~${\cal Q}$.
We use the standard abbreviations for the strong and weak operator topologies, i.e.~SOT and WOT respectively.

\section{Compact quantum groups}\label{Section2}

\subsection{Definition and representation of compact quantum groups}\label{2.1}

A~compact quantum group (CQG for short) is a~unital $C^{\ast}$-algebra ${\cal Q}$ with a~coassociative coproduct
(see~\cite{Van,CQG})~$\Delta$ from ${\cal Q}$ to ${\cal Q} \,\hat{\otimes}\, {\cal Q}$ such that each of the linear spans of
$\Delta({\cal Q})({\cal Q}\otimes 1)$ and that of $\Delta({\cal Q})(1\otimes {\cal Q})$ is norm-dense in ${\cal Q}
\,\hat{\otimes}\, {\cal Q}$.
From this condition, one can obtain a~canonical dense unital $\ast$-subalgebra ${\cal Q}_0$ of ${\cal Q}$ on which
linear maps~$\kappa$ and~$\epsilon$ (called the antipode and the counit respectively) are def\/ined making the above
subalgebra a~Hopf $\ast$-algebra.
In fact, we shall always choose this dense Hopf $\ast$-algebra to be the algebra generated by the `matrix coef\/f\/icients'
of the (f\/inite-dimensional) irreducible non-degenerate representations (to be def\/ined shortly) of the CQG.
The antipode is an anti-homomorphism and also satisf\/ies $\kappa(a^*)=(\kappa^{-1}(a))^*$ for $a \in {\cal Q}_{0}$.
A~morphism from a~CQG ${\cal Q}_1$ to another CQG ${\cal Q}_2$ is a~unital $\ast$-homomorphism~$\pi$ such that $\Delta_2
\circ \pi=(\pi \otimes \pi) \circ \Delta_1$, where $\Delta_i$ denotes the coproduct of ${\cal Q}_i$ ($i=1,2$).
If~$\pi$ is surjective, we say that ${\cal Q}_2$ is a~quantum subgroup of ${\cal Q}_1$ and write ${\cal Q}_{2}<{\cal
Q}_{1}$.
We shall use standard Swedler notation for Hopf algebras, i.e.~we shall write $\Delta(a)=a_{(1)}\otimes a_{(2)}$ and
$(\Delta\otimes {\rm id})\Delta(a)=a_{(1)}\otimes a_{(2)}\otimes a_{(3)}$ for $a\in{\cal Q}_{0}$.

It is known that there is a~unique state~$h$ on a~CQG ${\cal Q}$ (called the Haar state) which is bi-invariant in the
sense that $({\rm id} \otimes h)\circ \Delta(a)=(h \otimes {\rm id}) \circ \Delta(a)=h(a)1$ for all $a\in{\cal Q}$.
The Haar state need not be faithful in general, though it is always faithful on ${\cal Q}_0$ at least.
Given the Hopf $\ast$-algebra ${\cal Q}_{0}$, there can be several CQG's which have this $\ast$-algebra as the Hopf
$\ast$-algebra generated by the matrix elements of f\/inite-dimensional representations.
However there exists a~largest such CQG ${\cal Q}^{u}$, called the universal CQG corresponding to ${\cal Q}_{0}$.
It is obtained as the universal enveloping $C^{\ast}$-algebra of ${\cal Q}_{0}$.
We also say that a~CQG ${\cal Q}$ is universal if ${\cal Q}={\cal Q}^{u}$.
For details the reader is referred to~\cite{locallycompact}.
The $C^*$-completion ${\cal Q}_r$ of ${\cal Q}_0$ in the norm of ${\cal B}(L^2(h))$ (GNS space associated to~$h$) is
a~CQG and called the reduced quantum group corresponding to ${\cal Q}$.
If~$h$ is faithful then ${\cal Q}$ and ${\cal Q}_r$ coincide.
In general there will be a~surjective CQG morphism from ${\cal Q}$ to ${\cal Q}_r$ identifying the latter as a~quantum
subgroup of the former.

There is also a~von Neumann algebraic framework of quantum groups suitable for development of the theory of
locally compact quantum groups (see~\cite{kustermans,vaes} and references therein).
In this theory, the von Neumann algebraic version of CQG is a~von Neumann algebra ${\cal M}$ with a~coassociative,
normal, injective coproduct map~$\Delta$ from ${\cal M}$ to ${\cal M} \bar{\otimes}{\cal M}$ and a~faithful, normal,
bi-invariant state~$\psi$ on~${\cal M}$.
Indeed, given a~CQG~${\cal Q}$, the weak closure ${\cal Q}_r^{\prime \prime}$ of the reduced quantum group in the GNS
space of the Haar state is a~von Neumann algebraic compact quantum group.

Let ${\cal H}$ be a~Hilbert space.
Consider the multiplier algebra ${\cal M}({\cal K}({\cal H})\,\hat{\otimes}\, {\cal Q})$ which can also be identif\/ied with
${\cal L}({\cal H}\otimes{\cal Q})$, the $C^{\ast}$-algebra of ${\cal Q}$ linear, adjointable maps on ${\cal
H}\otimes{\cal Q}$ (see~\cite{Lance}).
This algebra has two natural embeddings into ${\cal M}({\cal K}({\cal H})\,\hat{\otimes}\, {\cal Q}\,\hat{\otimes}\, {\cal Q})$.
The f\/irst one is obtained by extending the map $x\mapsto x\otimes 1$.
The second one is obtained by composing this map with the f\/lip on the last two factors.
We will write $w^{12}$ and $w^{13}$ for the images of an element $w\in {\cal M}({\cal K}({\cal H})\,\hat{\otimes}\, {\cal
Q})$ by these two maps respectively.
Note that if ${\cal H}$ is f\/inite-dimensional then ${\cal M}({\cal K}({\cal H})\,\hat{\otimes}\, {\cal Q})$ is isomorphic to
${\cal B}({\cal H})\otimes {\cal Q}$ (we don't need any topological completion).

\begin{dfn}
Let $({\cal Q},\Delta)$ be a~CQG.
A~unitary representation~$U$ of ${\cal Q}$ on a~Hilbert space ${\cal H}$ is a~unitary element of ${\cal M}({\cal
K}({\cal H})\,\hat{\otimes}\,{\cal Q})$ such that $({\rm id}\otimes \Delta)U= U^{12}U^{13}$.
\end{dfn}

Identifying ${\cal M}({\cal K}({\cal H})\,\hat{\otimes}\,{\cal Q})$ with ${\cal L}({\cal H}\otimes{\cal Q})$, we can also
view~$U$ as a~${\cal Q}$ linear map on ${\cal H}\otimes{\cal Q}$.
Clearly~$U$ is determined by its restriction on the subspace ${\cal H}\otimes 1$ and by a~slight abuse of notation we
shall denote the restriction also by~$U$, which is a~$\mathbb{C}$-linear map $\xi\mapsto U(\xi\otimes q)\in{\cal
H}\otimes{\cal Q}$ for $\xi\in{\cal H}$, $q\in{\cal Q}$.
This gives an equivalent picture of a~unitary representation.
This is easy to see that the restriction $U:{\cal H}\rightarrow{\cal H}\otimes{\cal Q}$ satisf\/ies
\begin{enumerate}\itemsep=0pt
\item[1)] $\langle\langle U(\xi),U(\eta)\rangle \rangle =\langle \xi,\eta\rangle 1_{{\cal Q}}$ for $\xi,\eta\in{\cal H}$,

\item[2)] $(U\otimes {\rm id})U=({\rm id}\otimes\Delta)U$,

\item[3)]
$\operatorname{Sp} \{U(\xi)q:\xi\in{\cal H},q\in{\cal Q}\}$ is dense in ${\cal H}\otimes{\cal Q}$.
\end{enumerate}

\begin{dfn}
A~closed subspace ${\cal H}_{1}$ of ${\cal H}$ is said to be invariant if $U({\cal H}_{1})$ $(\equiv U({\cal H}_{1}\otimes
1_{{\cal Q}}))\subset {\cal H}_{1}\otimes{\cal Q}$.
Equivalently, the orthogonal projection~$p$ on ${\cal H}_{1}$ satisf\/ies $U(p\otimes 1)=(p\otimes 1)U$.
A~unitary representation~$U$ of a~CQG is said to be irreducible if there is no proper invariant subspace.
\end{dfn}

It is a~well known fact that every irreducible unitary representation is f\/inite-dimensional.

We denote by $\operatorname{Rep}({\cal Q})$ the set of inequivalent irreducible unitary representations of ${\cal Q}$.
For $\pi\in \operatorname{Rep}({\cal Q})$, let $d_{\pi}$ and $\{q^{\pi}_{jk}: j,k=1,\dots,d_{\pi} \}$ be the dimension and matrix
coef\/f\/icients of the corresponding f\/inite-dimensional representation, say $U_{\pi}$ respectively.
For each $\pi\in \operatorname{Rep}({\cal Q})$, we have a~unique $d_{\pi}\times d_{\pi}$ complex matrix $F_{\pi}$ such that
\begin{enumerate}\itemsep=0pt
\item[(i)] $F_{\pi}$ is positive and invertible with $\operatorname{Tr}(F_{\pi})=\operatorname{Tr}(F_{\pi}^{-1})=M_{\pi}>0$ (say),

\item[(ii)] $h(q_{ij}^{\pi}q_{kl}^{\pi^{\ast}})= \frac{1}{M_{\pi}}\delta_{ik}F_{\pi}(j,l)$.
\end{enumerate}

Corresponding to $\pi\in \operatorname{Rep}({\cal Q})$, let $\rho^{\pi}_{sm}$ be the linear functional on ${\cal Q}$ given~by
$\rho^{\pi}_{sm}(x)=h(x^{\pi}_{sm}x)$, $s,m=1,\dots ,d_{\pi}$ for $x\in {\cal Q}$,
where $x^{\pi}_{sm}=(M_{\pi})q^{\pi *}_{km}(F_{\pi}(k,s))$.
Also let $\rho^{\pi}=\sum\limits_{s=1}^{d_{\pi}}\rho^{\pi}_{ss}$.
Given a~unitary representation~$V$ on a~Hilbert space ${\cal H}$, we get a~decomposition of ${\cal H}$ as
\begin{gather*}
{\cal H}=\oplus_{\pi\in \operatorname{Rep}({\cal Q}),1\leq i\leq m_{\pi}}{\cal H}^{\pi}_{i},
\end{gather*}
where $m_{\pi}$ is the multiplicity of~$\pi$ in the representation~$V$ and $V|_{{\cal H}^{\pi}_{i}}$ is same as the
representa\-tion~$U_{\pi}$.
The subspace ${\cal H}^{\pi}=\oplus_{i}{\cal H}^{\pi}_{i}$ is called the spectral subspace of type~$\pi$ corresponding
to the irreducible representation~$\pi$.
It is nothing but the image of the spectral projection given by $({\rm id}\otimes \rho^{\pi})V$.

Recall from~\cite{Tomita}, the modular operator $\Phi=S^{\ast}S$, where~$S$ is the anti unitary acting on the
$L^{2}(h)$ (where $L^{2}(h)$ is the GNS space of ${\cal Q}$ corresponding to the Haar state on which ${\cal Q}$ has left
regular representation) given by $S(a.1):=a^{\ast}.1$ for $a\in{\cal Q}$.
The one parameter modular automorphism group (see~\cite{Tomita}) say $\Theta_{t}$, corresponding to the state~$h$ is
given by $\Theta_{t}(a)= \Phi^{it}a\Phi^{-it}$.
Note that here we have used the symbol~$\Phi$ for the modular operator as~$\Delta$ has been used for the coproduct.
From~(ii), we see that
\begin{gather}
\Phi|_{L^{2}(h)^{\pi}_{i}}=F_{\pi} \quad \text{for~all} \ \pi \ \text{and} \ i.
\label{eq1}
\end{gather}
In particular~$\Phi$ maps $L^{2}(h)^{\pi}_{i}$ into $L^{2}(h)^{\pi}_{i}$ for all~$i$.

\subsection{Action on von Neumann algebras by conjugation\\ of unitary representation}

We now discuss an analogue of `action' (as in~\cite{Podles}) in the context of von Neumann algebra implemented~by
a~unitary representation of the CQG.
Given a~unitary representation~$V$ of a~CQG~${\cal Q}$ on a~Hilbert space ${\cal H}$, often we consider the $\ast$
homomorphism ${\rm ad}_{V}$ on ${\cal B}({\cal H})$ or on some suitable von Neumann subalgebra ${\cal M}$ of it.
We say ${\rm ad}_{V}$ leaves ${\cal M}$ invariant if $({\rm id}\otimes \phi){\rm ad}_{V}({\cal M})\subset {\cal M}$ for every
state~$\phi$ of ${\cal Q}$.
Then taking $\rho^{\pi}$ as above, we def\/ine ${\cal M}^{\pi}=P_{\pi}({\cal M})$, where $P_{\pi}=({\rm id}\otimes
\rho^{\pi})\circ{\rm ad}_{V}:{\cal M}\rightarrow{\cal M}$ is the spectral projection corresponding to the
representation~$\pi$.
Clearly each $P_{\pi}$ is SOT continuous.
We def\/ine ${\cal M}_{0}:=\operatorname{Sp}\{{\cal M}^{\pi};\pi\in \operatorname{Rep}({\cal Q})\}$, which is called the spectral subalgebra.
Then we have the following:

\begin{ppsn}\label{SOT}
${\cal M}_{0}$ is dense in ${\cal M}$ in any of the natural locally convex topologies of ${\cal M}$, i.e.~${\cal M}_0^{\prime \prime}={\cal M}$.
\end{ppsn}

\begin{proof}
This result must be quite well-known and available in the literature but we could not f\/ind it written in this form, so
we give a~very brief sketch.
The proof is basically almost verbatim adaptation of some arguments in~\cite{kustermans} and~\cite{vaes}.
First, observe that the spectral algebra~${\cal M}_0$ remains unchanged if we replace~${\cal Q}$ by the reduced quantum
group ${\cal Q}_r$ which has the same irreducible representations and dense Hopf $\ast$-algebra~${\cal Q}_0$.
This means we can assume without loss of generality that the Haar state is faithful.
The injective normal map $\beta:={\rm ad}_{V}$ restricted to ${\cal M}$ can be thought of as an action of the quantum
group (in the von Neumann algebra setting as in~\cite{kustermans,vaes}) ${\cal Q}_r^{\prime \prime}$ (where
the double commutant is taken in the GNS space of the Haar state) and it follows from the results in~\cite{vaes} about
the implementability of locally compact quantum group actions that there is a~faithful normal state, say~$\phi$, on
${\cal M}$ such that~$\beta$ is ${\cal Q}_r^{\prime \prime}$-invariant, i.e.~$(\phi \otimes {\rm id})(\beta(a))=\phi(a)1$ $\forall\, a \in {\cal M}$.
We can replace ${\cal M}$, originally imbedded in ${\cal B}({\cal H})$, by its isomorphic image (to be denoted by ${\cal
M}$ again) in ${\cal B}(L^2(\phi))$ and as the ultra-weak topology is intrinsic to a~von Neumann algebra, it will
suf\/f\/ice to argue the ultra-weak density of~${\cal M}_0$ in~${\cal M} \subset {\cal B}(L^2(\phi))$.
To this end, note that Vaes has shown in~\cite{vaes} that $\beta: {\cal M} \rightarrow {\cal M} \otimes {\cal
Q}_r^{\prime \prime}$ extends to a~unitary representation on~$L^2(\phi)$, which implies in particular that~${\cal M}_0$
is dense in the Hilbert space~$L^2(\phi)$.
From this the ultra-weak density follows by standard arguments very similar to those used in the proof of Proposition~1.5 of~\cite{kustermans}, applying Takesaki's theorem about existence of conditional expectation.
For the sake of completeness let us sketch it brief\/ly.
Using the notations of~\cite{vaes} and noting that $\delta=1$ for a~CQG, we get from Proposition~2.4 of~\cite{vaes} that~$V_{\phi}$ commutes with the positive self adjoint operator $\nabla_{\phi}\otimes Q$ where $\nabla_{\phi}$ denotes the
modular operator, i.e.~generator of the modular automorphism group $\sigma^{\phi}_{t}$ of the normal state~$\phi$.
Clearly, this implies that $\beta:={\rm ad}_{V_{\phi}}$ satisf\/ies the following:
\begin{gather*}
\beta\circ \sigma^{\phi}_{t}=\sigma^{\phi}_{t}\otimes \tau_{-t},
\end{gather*}
where $\tau_{t}$ is the automorphism group generated by $Q^{-1}$.
Next, as in Proposition~1.5 of~\cite{kustermans}, consider the ultra-strong $\ast$ closure ${\cal M}_{l}$ of the
subspace spanned by the elements of the form $({\rm id}\otimes\omega)(\beta(x))$, $x\in M$, $\omega$~is a~bounded normal
functional on ${\cal Q}_{r}^{\prime\prime}$.
It is enough to prove that ${\cal M}_{l}={\cal M}$, as that will prove the ultra-strong $\ast$ density (now also the
ultra-weak density) of ${\cal M}_{l}$ in ${\cal M}$.
This is clearly a~von Neumann subalgebra as~$\beta$ is coassociative, and
$\sigma^{\phi}_{t}({\rm id}\otimes\omega)(\beta(x))=({\rm id}\otimes \omega\circ \tau_{t})(\beta(\sigma^{\phi}_{t}(x)))$.
Then by Takesaki's theorem \cite[Theorem~10.1]{Takesaki} there exists a~unique normal faithful conditional expectation~$E$
from ${\cal M}$ to ${\cal M}_{l}$ satisfying $E(x)P=PxP$ where~$P$ is the orthogonal projection as in~\cite{kustermans}.
Clearly, the range of~$P$ contains elements of the form $({\rm id}\otimes\omega)(\beta(x))$.
So in particular it contains ${\cal M}_{0}$, which is dense in $L^{2}(\phi)$.
Thus $P=1$ and $E(x)=x$ proving ${\cal M}_{l}={\cal M}$.
\end{proof}

From this, it is easy to conclude the following by verbatim adaptation of the arguments of~\cite{Podles}
and~\cite{Soltan_action}.
\begin{lmma}\label{new1}\quad\sloppy
\begin{enumerate}\itemsep=0pt
\item[$1.$]
${\rm ad}_{V}|_{{\cal M}_{0}}$ is algebraic, i.e.~${\rm ad}_{V}({\cal M}_{0})\subset {\cal M}_{0}\otimes_{\rm alg} {\cal Q}_{0}$.

\item[$2.$]
${\cal M}_{0}$ is the maximal subspace over which ${\rm ad}_{V}$ is algebraic,
i.e.~${\cal M}_{0}=\{x\in{\cal M}\,|\, {\rm ad}_{V}(x)\in{\cal M}\otimes_{\rm alg}{\cal Q}_{0}\}$.

\item[$3.$]
If ${\cal M}_{1}\subset {\cal M}$ is SOT dense $\ast$-subalgebra such that ${\rm ad}_{V}$ leaves ${\cal M}_{1}$
invariant, then $\operatorname{Sp}\{P_{\pi}({\cal M}_{1})\,| \, \pi\in \operatorname{Rep}({\cal Q})\}$ is SOT dense in ${\cal M}_{0}$.
\end{enumerate}
\end{lmma}

\begin{proof}
(1) Recall the linear functionals $\rho^{\pi}_{sm}$ from Section~\ref{2.1}.
Let $P^{\pi}_{sm}:=({\rm id}\otimes\rho^{\pi}_{sm})\circ \Delta:{\cal Q}\rightarrow{\cal Q}$ and $E^{\pi}_{sm}:=({\rm id}\otimes
\rho^{\pi}_{sm}){\rm ad}_{V}:{\cal M}\rightarrow{\cal M}$ for $\pi\in \operatorname{Rep}({\cal Q})$ and $s,m=1,\dots ,d_{\pi}$.
Then $P^{\pi}_{sm}({\cal Q})\subset \operatorname{Sp} \{q^{\pi}_{is}:i=1,\dots ,d_{\pi}\}$.
It follows from the arguments of~\cite{Podles},
\begin{gather}\label{eq2}
E^{\pi}_{sm}E^{\pi^{\prime}}_{ij}=\delta_{mi}\delta_{\pi\pi^{\prime}}E^{\pi}_{is}
\end{gather}
for $\pi,\pi^{\prime}\in \operatorname{Rep}({\cal Q})$.
Let ${\cal M}_{\pi s}=\sum\limits_{s=1}^{d_{\pi}}E^{\pi}_{ss}({\cal M})$, we have by~\eqref{eq2},
\begin{gather}\label{eq3}
{\cal M}^{\pi}=\oplus_{s=1}^{d_{\pi}}{\cal M}_{\pi s}.
\end{gather}

Let $\{e_{\pi i1}\}_{i\in {\cal J}}$ be a~basis for ${\cal M}_{\pi 1}$ for $\pi\in \operatorname{Rep}({\cal Q})$, and $e_{\pi is}:=
E^{\pi}_{s1}(e_{\pi i1})$, $s=1,\dots ,d_{\pi}$, $i\in {\cal J}$.
Then by~\eqref{eq2}, $\{e_{\pi is}\}_{i\in {\cal J}}$ is a~basis for ${\cal M}_{\pi s}$.
Letting ${\cal M}^{\pi}_{i}= \operatorname{Sp}\{e_{\pi is}:s=1,\dots ,d_{\pi}\}$ and using~\eqref{eq3}, we get
${\cal M}^{\pi}=\oplus_{i\in {\cal J}}{\cal M}^{\pi}_{i}$.

Now
\begin{gather*}
\begin{split}
& {\rm ad}_{V}(e_{\pi is}) =  ({\rm id}\otimes({\rm id}\otimes\rho^{\pi}_{ss})\Delta){\rm ad}_{V}(e_{\pi is})\\
& \hphantom{{\rm ad}_{V}(e_{\pi is})}{}
 =  ({\rm id}\otimes P^{\pi}_{ss}){\rm ad}_{V}(e_{\pi is})
 \subset  {\cal M}\otimes \operatorname{Sp} \{q^{\pi}_{js};j=1,\dots ,d_{\pi}\}.
\end{split}
\end{gather*}
Hence ${\rm ad}_{V}(e_{\pi is})=\sum\limits_{j=1}^{d_{\pi}} x_{(\pi)i(s)j}\otimes q^{\pi}_{js}$ for some
$x_{(\pi)i(s)j}\in {\cal M}$.
Now applying $({\rm id}\otimes \rho^{\pi}_{ks})$ for $k=1,\dots ,d_{\pi}$ on both sides, we get $x_{(\pi)i(s)k}=e_{\pi ik}$,
i.e.~${\rm ad}_{V}(e_{\pi is})=\sum\limits_{j=1}^{d_{\pi}}e_{\pi ij}\otimes q^{\pi}_{js}$
proving ${\rm ad}_{V}({\cal M}_{0})\subset {\cal M}_{0}\otimes_{\rm alg}{\cal Q}_{0}$.

(2)~Here we reproduce the arguments of~\cite{Soltan_action} for the sake of completeness.
For any $b\in{\cal M}$ such that ${\rm ad}_{V}(b)\subset{\cal M}\otimes_{\rm alg}{\cal Q}_{0}$, we have
\begin{gather*}
{\rm ad}_{V}(b)=\sum\limits_{\pi\in {\cal S}\subset \operatorname{Rep}({\cal Q})}\sum\limits_{i,j=1}^{d_{\pi}}b^{\pi}_{ij}\otimes
q^{\pi}_{ij},
\end{gather*}
where ${\cal S}$ is a~f\/inite subset of $\operatorname{Rep}({\cal Q})$.
Observe that for each $\pi\in{\cal S}$, $i,j=1,\dots ,d_{\pi}$, $b^{\pi}_{ij}=({\rm id}\otimes \rho^{\pi}_{ij})\circ {\rm
ad}_{V}(b)\in{\cal M}_{0}$.
As $({\rm ad}_{V}\otimes {\rm id})\circ{\rm ad}_{V}(b)=({\rm id}\otimes \Delta)\circ{\rm ad}_{V}(b)$, we have
\begin{gather*}
\sum\limits_{\pi\in{\cal S}\subset \operatorname{Rep}({\cal Q})}\sum\limits_{i,j=1}^{d_{\pi}}{\rm ad}_{V}(b^{\pi}_{ij})\otimes
q^{\pi}_{ij}=\sum\limits_{\pi\in{\cal S}\subset \operatorname{Rep}({\cal Q})}\sum\limits_{i,j,s=1}^{d_{\pi}}b^{\pi}_{ij}\otimes
q^{\pi}_{is}\otimes q^{\pi}_{sj}.
\end{gather*}
Applying $({\rm id}\otimes {\rm id}\otimes \rho^{\pi}_{kl})$ on both sides, we get
\begin{gather}\label{eq4}
{\rm ad}_{V}(b^{\pi}_{kl})=\sum\limits_{i=1}^{d_{\pi}} b^{\pi}_{il}\otimes q^{\pi}_{ik}.
\end{gather}
Now if we take $b^{\prime}=\sum\limits_{\pi\in{\cal S}}\sum\limits_{i=1}^{d_{\pi}}b^{\pi}_{ii}$, we get by~\eqref{eq4},
$
{\rm ad}_{V}(b^{\prime})=\sum\limits_{\pi\in{\cal S}}\sum\limits_{k=1}^{d_{\pi}}b^{\pi}_{ki}\otimes q^{\pi}_{ki}$.
So ${\rm ad}_{V}(b)={\rm ad}_{V}(b^{\prime})$ and as ${\rm ad}_{V}$ is one-one, $b=b^{\prime}\in{\cal M}_{0}$.

(3) follows from (2) and the SOT continuity of each $P_{\pi}$.
\end{proof}

For $x\in{\cal M}_{0}$, we shall use the natural analogue of Swedler's notation, i.e.~write
${\rm ad}_{V}(x)=x_{(0)}\otimes x_{(1)}$.

\subsection{Dual of a~compact quantum group}\label{2.3}

Let $({\cal Q},\Delta)$ be a~compact quantum group and ${\cal Q}_{0}$ be its canonical dense Hopf $\ast$-algebra as in
Section~\ref{2.1}.
We def\/ine $\hat{{\cal Q}_{0}}$ to be the space of linear functionals on ${\cal Q}$ def\/ined by $q\rightarrow h(aq)$, for
$a\in {\cal Q}_{0}$, where~$h$ is the Haar state on ${\cal Q}$.
Then $\hat{{\cal Q}_{0}}$ is a~subspace of the dual ${\cal Q}^{\prime}$.
For $\pi\in \operatorname{Rep}({\cal Q})$, let $F_{\pi}$ be the positive invertible matrix as in Section~\ref{2.1}.
Recall also the linear functionals $\rho^{\pi}_{sm}$'s from Section~\ref{2.1}.
Then we have from Lemma~8.1 of~\cite{Van}, $\rho^{\pi}_{pr}(q^{\pi}_{pr})=1$ and $\rho^{\pi}_{pr}$ is zero on all other
matrix elements.

Note that (see~\cite{Van}) $\hat{{\cal Q}_{0}}$ is a~$\ast$-subalgebra of ${\cal Q}^{\prime}$ and as a~$\ast$-algebra
it is nothing but the algebraic direct sum of matrix algebras of the form $\oplus_{\pi\in \operatorname{Rep}({\cal Q})}M_{d_{\pi}}$
where $M_{d_{\pi}}$ is the full matrix algebra of size $d_{\pi}\times d_{\pi}$ and the matrix unit $m^{\pi}_{pq}$
identif\/ied with $\rho^{\pi}_{pq}$ is given by $|e^{\pi}_{p}\rangle\langle e^{\pi}_{q}|$, where
$\{e^{\pi}_{1},\dots ,e^{\pi}_{d_{\pi}}\}$ is the standard orthonormal basis.
Clearly we get a~non-degenerate pairing between $\hat{{\cal Q}_{0}}=\oplus_{\pi\in \operatorname{Rep}({\cal Q})} M_{d_{\pi}}$ and
${\cal Q}_{0}$ given by $\langle q^{\pi}_{ij},m^{\pi^{\prime}}_{ps}\rangle=\delta_{ip}\delta_{js}\delta_{\pi\pi^{\prime}}$.
In fact it can be shown that $\hat{{\cal Q}_{0}}$ is a~multiplier Hopf $\ast$-algebra in the sense of~\cite{mult}.
It has a~unique $C^{\ast}$ norm and by taking the completion we get a~$C^{\ast}$-algebra.
It is called the dual discrete quantum group of ${\cal Q}$ denoted by $\hat{{\cal Q}}$.
It has a~coassociative $\ast$-homomorphism $\hat{\Delta}:\hat{{\cal Q}}\rightarrow {\cal M}(\hat{{\cal
Q}}\,\hat{\otimes}\,\hat{{\cal Q}})$ called the dual coproduct satisfying
$({\rm id}\otimes\hat{\Delta})\hat{\Delta}=(\hat{\Delta}\otimes {\rm id})\hat{\Delta}$.

\begin{thm}
\label{mult_rep}
Let~$U$ be a~unitary representation of a~CQG ${\cal Q}$ on a~Hilbert space ${\cal H}$.
Then $\Pi_{U}:\hat{{\cal Q}}\rightarrow {\cal B}({\cal H})$ defined by $\Pi_{U}(\omega)(\xi):=({\rm id}\otimes \omega)U(\xi)$
is a~non degenerate $\ast$-homomorphism and hence extends as a~$\ast$-homomorphism from ${\cal M}(\hat{{\cal Q}})$ to
${\cal B}({\cal H})$.
\end{thm}

\begin{proof}
Consider the spectral decomposition ${\cal H}=\oplus_{\pi\in {\cal I},1\leq i\leq m_{\pi}}{\cal H}^{\pi}_{i}$,
$U|_{{\cal H}^{\pi}_{i}}$, $i=1,\dots ,m_{\pi}$ is equivalent to the irreducible representation of type~$\pi$.
Moreover f\/ix orthonormal basis $e^{\pi}_{ij}$, $j=1,\dots ,d_{\pi}$, $i=1,\dots ,m_{\pi}$ for ${\cal H}^{\pi}_{i}$ such that
$
U(e^{\pi}_{ij})=\sum\limits_{k}e^{\pi}_{ik}\otimes q^{\pi}_{kj}
$
for all $\pi\in \operatorname{Rep}({\cal Q})$.
Now for a~f\/ixed $\pi\in \operatorname{Rep}({\cal Q})$, $p,r=1,\dots ,d_{\pi}$ observe that $\Pi_{U}(\rho^{\pi}_{pr})(\xi)=0$ for all
$\xi\in{\cal H}^{\pi^{\prime}}_{i}$ and for $\pi\neq\pi^{\prime}$.
Also $\Pi_{U}(\rho^{\pi}_{pr})(e^{\pi}_{ij})=\delta_{jr}e^{\pi}_{ip}$, i.e.~$\Pi_{U}(\rho^{\pi}_{pr})|_{{\cal H}^{\pi}_{i}}$
is nothing but the rank one operator $|e^{\pi}_{ip}\rangle\langle e^{\pi}_{ir}|$.
This proves that $\Pi_{U}(\omega)$ is bounded for $\omega\in\hat{{\cal Q}_{0}}$, and moreover identifying $\hat{{\cal
Q}_{0}}$ with the direct sum of matrix algebras $\oplus_{\pi\in \operatorname{Rep}({\cal Q})}M_{d_{\pi}}$, we see that $\Pi_{U}$ is
nothing but the map which sends $X\in M_{d_{\pi}}$ to $X\otimes 1_{\mathbb{C}^{m\pi}}$ in ${\cal B}({\cal H})$.
This proves that $\Pi_{U}$ extends to a~non-degenerate $\ast$-homomorphism.
\end{proof}

\section{Quantum isometry groups}\label{Section3}

Before def\/ining quantum isometry group, we recall the def\/inition of a~spectral triple of compact type.
\begin{dfn}
A~triple $({\cal A}^{\infty},{\cal H},{\cal D})$, where ${\cal H}$ is Hilbert space, ${\cal A}^{\infty}$ is a~unital
$\ast$-subalgebra of~${\cal B}({\cal H})$ and~${\cal D}$ is an unbounded operator on~${\cal H}$ is called a~spectral
triple of compact type if
\begin{enumerate}\itemsep=0pt
\item[(i)] $[{\cal D},a]\in{\cal B}({\cal H})$ for all $a\in {\cal A}^{\infty}$,
\item[(ii)] ${\cal D}$ has compact resolvents.
\end{enumerate}
\end{dfn}

Now let us brief\/ly sketch the formulation of quantum isometry groups of spectral triples given by Goswami and Bhowmick
in~\cite{qorient,Goswami}.
\begin{dfn}
Let $({\cal A}^{\infty},{\cal H},{\cal D})$ be a~spectral triple of compact type (a la Connes).
Consider the category ${\bf Q}^{\prime}({\cal D})\equiv {\bf Q}^{\prime}({\cal A}^{\infty},{\cal H},{\cal D})$ whose objects are
the triples $({\cal Q},\Delta, U)$ where $({\cal Q},\Delta)$ is a~CQG having a~unitary representation~$U$ on the Hilbert
space ${\cal H}$ such that~$U$ commutes with $({\cal D} \otimes 1_{{\cal Q}})$.
A~morphism from $({\cal Q},\Delta, U)$ to $({\cal Q}^{\prime},\Delta^{\prime}, U^{\prime})$ is a~CQG morphism $\psi:{\cal Q}\rightarrow
{\cal Q}^{\prime}$ such that $U^{\prime}=({\rm id}\otimes \psi)U$.
\end{dfn}
\begin{dfn}
For a~positive (possibly unbounded) operator~$R$ on ${\cal H}$ which commutes with~${\cal D}$, consider the subcategory
${\bf Q}^{\prime}_{R}({\cal D})\equiv {\bf Q}^{\prime}_{R}({\cal A}^{\infty},{\cal H},{\cal D})$ whose objects are triples $({\cal
Q},\Delta,U)$ as before with the additional requirement that ${\rm ad}_{U}$ preserves the functional $\tau_{R}$ given~by
$\tau_{R}(x)=\operatorname{Tr}(Rx) (x\in {\cal E}_{{\cal D}})$, in the sense that $(\tau_{R}\otimes {\rm id}){\rm
ad}_{U}(x)=\tau_{R}(x).1_{{\cal Q}}$ for $x\in{\cal E}_{{\cal D}}$, where ${\cal E}_{{\cal D}}$ is the weakly dense
$\ast$-subalgebra of ${\cal B}({\cal H})$ spanned by the rank one operators of the form $|\xi\rangle\langle \eta|$ where~$\xi$
and~$\eta$ are eigenvectors of~${\cal D}$.
\end{dfn}
\begin{rmrk}
We refer the functional $\tau_{R}$ as the~$R$-twisted volume and say that any object of~$Q^{\prime}_{R}({\cal D})$ preserves
the~$R$-twisted volume.
\end{rmrk}

\begin{ppsn}[\protect{\cite[Theorem 2.14]{qorient}}] For a~spectral triple of compact type, the universal object
in the category ${\bf Q}^{\prime}_{R}({\cal D})$ exists and is denoted~by
$\widetilde{{\rm QISO}^{+}_{R}}({\cal A}^{\infty},{\cal H},{\cal D})$.
\end{ppsn}

\begin{dfn}
The largest Woronowicz subalgebra of $\widetilde{{\rm QISO}^{+}_{R}}({\cal A}^{\infty},{\cal H},{\cal D})$ on which ${\rm
ad}_{U}$ is faithful ($U$~is the corresponding unitary representation) is called the quantum group of orientation
preserving Riemannian isometry of the~$R$-twisted spectral triple and denoted by ${\rm QISO}^{+}_{R}({\cal A}^{\infty}$, ${\cal
H},{\cal D})$ or simply by ${\rm QISO}^{+}_{R}({\cal D})$.
\end{dfn}

\begin{rmrk}\label{remark3.7}
${\bf Q}^{\prime}({\cal D})$ may not have a~universal object in general, but if it exists we shall denote it~by
$\widetilde{{\rm QISO}^{+}}({\cal A}^{\infty},{\cal H},{\cal D})$ and the corresponding largest Woronowicz subalgebra~by
${\rm QISO}^{+}({\cal A}^{\infty}$, ${\cal H},{\cal D})$ or simply $Q{\rm ISO}^{+}({\cal D})$.
\end{rmrk}

Let $({\cal Q},V)$ be an object in the category ${\bf Q}^{\prime}_{R}({\cal D})$.
We would like to give a~necessary and suf\/f\/icient condition on the unbounded operator~$R$ so that ${\rm ad}_{V}$
preserves the~$R$-twisted volume.
To this end decompose the Hilbert space ${\cal H}$, on which ${\cal D}$ acts, into f\/inite-dimensional eigenspaces of the
operator~${\cal D}$, i.e.~let ${\cal H}=\oplus_{k}{\cal H}_{k}$ where each~${\cal H}_{k}$ is a~f\/inite-dimensional
eigenspace for ${\cal D}$.
Since ${\cal D}$ commutes with~$V$, $V$~preserves each of the~${\cal H}_{k}$'s and on each~${\cal H}_{k}$, $V$~is
a~unitary representation of the compact quantum group~${\cal Q}$.
Then we have the decomposition of each~${\cal H}_{k}$ into the irreducibles, say
$
{\cal H}_{k}=\oplus_{\pi\in{\cal I}_{k}}\mathbb{C}^{d_{\pi}}\otimes\mathbb{C}^{m_{\pi,k}}$,
where $m_{\pi,k}$ is the multiplicity of the irreducible representation of type~$\pi$ on ${\cal H}_{k}$ and ${\cal
I}_{k}$ is some f\/inite subset of $\operatorname{Rep}({\cal Q})$.
Since~$R$ commutes with~$V$, $R$ preserves direct summands of~${\cal H}_{k}$.
Then we have the following:

\begin{thm}
\label{form_R}
${\rm ad}_{V}$ preserves the~$R$-twisted volume if and only if
$
R|_{{\cal H}_{k}}=\oplus_{\pi\in{\cal I}_{k}}F_{\pi}\otimes T_{\pi,k}$,
for some $T_{\pi,k}\in{\cal B}(\mathbb{C}^{m_{\pi,k}})$, where $F_{\pi}$'s are as in Section~{\rm \ref{2.1}}.
\end{thm}
\begin{proof}
Only if part:

let $\{e_{i}\}_{i=1}^{d_{\pi}}$ and $\{f_{j}\}_{j=1}^{m_{\pi,k}}$ be orthonormal bases for $\mathbb{C}^{d_{\pi}}$ and
$\mathbb{C}^{m_{\pi,k}}$ respectively.
Also let $R(e_{i}\otimes f_{j})=\sum\limits_{s,t}R(s,t,i,j)e_{s}\otimes f_{t}$.

We have $V^{\ast}(e_{i}\otimes f_{j}\otimes 1_{{\cal Q}})=\sum\limits_{k}e_{k}\otimes f_{j}\otimes
q_{ik}^{\ast}$.
We denote the restriction of the trace of ${\cal B}({\cal H})$ on $\mathbb{C}^{d_{\pi}}\otimes \mathbb{C}^{m_{\pi,k}}$
again by $\operatorname{Tr}$.
Let $a\in {\cal B}(\mathbb{C}^{d_{\pi}}\otimes \mathbb{C}^{m_{\pi,k}})$ and $\chi(a):=\operatorname{Tr}(a.R)$.
Then we have
\begin{gather*}
(\chi\otimes h){\rm ad}_{V}(a) =  \sum\limits_{i,j}\langle {V}^{\ast}(e_{i}\otimes f_{j}\otimes 1_{{\cal Q}}), (a\otimes
1){V}^{\ast}R(e_{i}\otimes f_{j})\rangle
\\
\hphantom{(\chi\otimes h){\rm ad}_{V}(a)}{}
 =  \sum\limits_{i,j,k,s,t,u}\langle e_{k}\otimes f_{j}\otimes q_{ik}^{\ast},R(s,t,i,j)a(e_{u}\otimes f_{t})\otimes
q_{su}^{\ast}\rangle
\\
\hphantom{(\chi\otimes h){\rm ad}_{V}(a)}{}
 =  \sum\limits_{i,j,k,s,t,u}\frac{R(s,t,i,j)}{M_{\pi}}\langle e_{k}\otimes f_{j},a(e_{u}\otimes f_{t})\rangle \delta_{is}F_{\pi}(k,u)
\\
\hphantom{(\chi\otimes h){\rm ad}_{V}(a)}{}
 =  \sum\limits_{i,j,k,t,u}\frac{R(i,t,i,j)}{M_{\pi}}\langle e_{k}\otimes f_{j},a(e_{u}\otimes f_{t})\rangle F_{\pi}(k,u).
\end{gather*}
On the other hand
\begin{gather*}
\chi(a) =    \operatorname{Tr}(a.R)
 =  \sum\limits_{i,j}\langle e_{i}\otimes f_{j},aR(e_{i}\otimes f_{j})\rangle
 =  \sum\limits_{k,j,u,t}R(u,t,k,j)\langle e_{k}\otimes f_{j},a(e_{u}\otimes f_{t})\rangle .
\end{gather*}
Thus $(\chi\otimes h){\rm ad}_{V}(a)=\chi(a)$ implies:
\begin{gather}\label{eq5}
  \sum\limits_{i,j,k,t,u}\frac{R(i,t,i,j)}{M_{\pi}}\langle e_{k}\otimes f_{j},a(e_{u}\otimes f_{t})\rangle F_{\pi}(k,u)
 = \sum\limits_{k,j,u,t}R(u,t,k,j)\langle e_{k}\otimes f_{j},a(e_{u}\otimes f_{t})\rangle .\!\!\!
\end{gather}
Now f\/ix $u_{0}$, $t_{0}$ and consider $a\in {\cal B}({\cal H})$ such that $a(e_{u_{0}}\otimes f_{t_{0}})=e_{p}\otimes
f_{q}$ and zero on the other basis elements.
Then from \eqref{eq5}, we get{\samepage
\begin{gather*}
   \sum\limits_{i,j,k}\frac{R(i,t_{0},i,j)}{M_{\pi}}\langle e_{k}\otimes f_{j},e_{p}\otimes f_{q}\rangle F_{\pi}(k,u_{0})
 = \sum\limits_{k,j}R(u_{0},t_{0},k,j)\langle e_{k}\otimes f_{j},e_{p}\otimes f_{q}\rangle,
\end{gather*}
which gives $\sum\limits_{i} \frac{R(i,t_{0},i,q)}{M_{\pi}}F_{\pi}(p,u_{0})=R(u_{0},t_{0},p,q)$.}

This proves that $R|_{{\cal H}_{k}}=\oplus_{\pi\in{\cal I}_{k}}F_{\pi}\otimes T_{\pi,k}$ with some $T_{\pi,k}\in{\cal
B}(\mathbb{C}^{m_{\pi,k}})$ given by $T_{\pi,k}(t_{0},q)=\sum\limits_{i}\frac{R(i,t_{0},i,q)}{M_{\pi}}$.

The if part is straightforward and was essentially done in~\cite{Goswami_entire_cyclic}.
\end{proof}

Let us note that in case of quantum group of orientation preserving Riemannian isometry,
$\widetilde{{\rm QISO}^{+}_{R}}({\cal D})$ really depends on the von Neumann algebra $({\cal A}^{\infty})^{\prime\prime}$, not on the
algebra ${\cal A}^{\infty}$ itself.
More precisely, we have the following proposition, which follows from the def\/inition of $\widetilde{{\rm QISO}^{+}_{R}}({\cal D})$ and needs no proof:
\begin{ppsn}
%\label{qiso_von}
Let $({\cal A}^{\infty}, {\cal H},{\cal D},R)$ be as in~{\rm \cite{qorient}}.
If we have a~SOT dense subalgebra ${\cal A}_{0}$ of $({\cal A}^{\infty})^{\prime\prime}\subset {\cal B}({\cal H})$ such that
$[{\cal D},a]\in {\cal B}({\cal H})$ for all $a\in {\cal A}_{0}$, then $({\cal A}_{0},{\cal H},{\cal D})$ is again
a~spectral triple and
$
\widetilde{{\rm QISO}^{+}_{R}}({\cal A}_{0},{\cal H},{\cal D})\cong \widetilde{{\rm QISO}^{+}_{R}}({\cal A}^{\infty}, {\cal H},{\cal D}),
$
and hence ${\rm QISO}^{+}_{R}({\cal A}_{0},{\cal H},{\cal D})\cong {\rm QISO}^{+}_{R}({\cal A}^{\infty}, {\cal H},{\cal D})$.
\end{ppsn}

\looseness=-1
In~\cite{Goswami}, for a~spectral triple satisfying some mild regularity conditions, a~noncommutative analogue
of the Hodge Laplacian was constructed which is denoted by ${\cal L}={\cal L}_{{\cal D}}$.
We denote the category whose objects are triples $({\cal S},\Delta,\alpha)$ where $({\cal S},\Delta)$ is a~CQG
and~$\alpha$ is an action of the CQG on the manifold commuting with the Laplacian (as in the sense of Def\/inition~2.11
of~\cite{Goswami}) by ${\bf Q}^{\prime}_{{\cal L}}$ and morphism between two such objects is a~CQG morphism intertwining the
actions (see~\cite{Goswami}).

It is shown in~\cite{Goswami} that under certain regularity conditions on the spectral triple, universal object
exists in this category and we denote the universal object by ${\rm QISO}^{{\cal L}}$ and call it the quantum isometry group
of the spectral triple.

\begin{rmrk}
Note that both $\widetilde{{\rm QISO}^{+}_{R}}({\cal D})$ and ${\rm QISO}^{{\cal L}}$ are universal CQG's in the sense of Section~\ref{2.1}.
\end{rmrk}

When the spectral triple is classical, i.e.~$(C^{\infty}(M),{\cal H},{\cal D})$ where~$M$ is a~compact,
connected, Riemannian spin manifold, it was shown in~\cite{rigidity} that the quantum isometry group coincides with the
classical isometry group.
More precisely we state the following results from~\cite[Corollary~11.10]{rigidity}:
\begin{ppsn}
\label{prevmain}
For a~classical spectral triple $(C^{\infty}(M),{\cal H},{\cal D})$,
$
{\rm QISO}^{{\cal L}}(C^{\infty}(M),{\cal H},{\cal D})\cong C({\rm ISO}(M))
$
and
$
{\rm QISO}^{+}_{I}(C^{\infty}(M),{\cal H},{\cal D})\cong C({\rm ISO}^{+}(M)).
$
\end{ppsn}

\section{Cocycle twisting}\label{Section4}

In Sections~\ref{4.1} and~\ref{4.2}, we recollect some well known facts about cocycle twist of CQG by dual unitary cocycles.
There is nothing new in these two subsections and we mostly state the results with appropriate references.
However, for convenience of the reader, we occasionally sketch a~few proofs and give some details.

\subsection{Cocycle twist of a~compact quantum group}\label{4.1}

Let ${\cal Q}$ be a~compact quantum group.
Recall from
Section~\ref{2.1} the dense Hopf $\ast$-algebra ${\cal Q}_{0}$ spanned by the matrix coef\/f\/icients of its
inequivalent irreducible representations.
Also recall the dual discrete quantum group $\hat{{\cal Q}}$ of ${\cal Q}$.
With these notations we have the following
\begin{dfn}
By a~dual unitary 2-cocycle~$\sigma$ of a~compact quantum group ${\cal Q}$, we mean a~unitary element of ${\cal
M}(\hat{{\cal Q}}\,\hat{\otimes}\,\hat{{\cal Q}})$ satisfying
\begin{gather*}
(1\otimes \sigma)({\rm id}\otimes \hat{\Delta})\sigma=(\sigma\otimes 1)(\hat{\Delta}\otimes {\rm id})\sigma.
\end{gather*}
\end{dfn}
A~dual unitary 2-cocycle~$\sigma$ is said to be normalized if $(1\otimes p_{\epsilon})\sigma=1\otimes p_{\epsilon}$ and
$(p_{\epsilon}\otimes 1)\sigma=p_{\epsilon}\otimes 1$, where $p_{\pi}$'s are the minimal projections of $\oplus_{\pi\in
\operatorname{Rep}({\cal Q})}M_{d_{\pi}}$.
Without loss of generality we can always assume a~cocycle to be normalized and hence forth in our paper all dual unitary
2-cocycles are normalized.
This means in particular $\sigma(a,1)=\sigma^{-1}(a,1)=\epsilon(a)$ for all $a\in{\cal Q}_{0}$.

Recall from Section~\ref{2.3}, the discrete quantum group $\hat{{\cal Q}_{0}}$ with the coproduct $\hat{\Delta}$
def\/ined by $\hat{\Delta}(\omega)(a\otimes b):= \omega(ab)$,
$w\in \hat{{\cal Q}_{0}}$  and  $a,b\in {\cal Q}_{0}$.
We can deform the discrete quantum group $\hat{{\cal Q}_{0}}$ using $\sigma\in M(\hat{{\cal Q}}\,\hat{\otimes}\,\hat{{\cal
Q}})$.
The $\ast$-algebraic structure do not change where the coproduct changes~by
\begin{gather*}
\hat{\Delta}_{\sigma}(\cdot)=\sigma\hat{\Delta}(\cdot)\sigma^{-1}.
\end{gather*}
Using the cocycle condition of~$\sigma$, it can be easily shown that $\hat{\Delta}_{\sigma}$ is again coassociative.
It can be shown that $(\hat{{\cal Q}_{0}}_{\sigma},\hat{\Delta}_{\sigma})$ is again a~discrete quantum group.
By Proposition~3.12 and discussions before Proposition~4.5 of~\cite{bichon}, we see that we have a~CQG dual to the
discrete quantum group $(\hat{{\cal Q}_{0}}_{\sigma},\hat{\Delta}_{\sigma})$.

Now~$\sigma$ viewed as a~linear functional on ${\cal Q}_{0}\otimes_{\rm alg}{\cal Q}_{0}$ satisf\/ies the cocycle
condition (see~\cite[p.~64]{twist})
\begin{gather*}
\sigma(b_{(1)},c_{(1)})\sigma (a,b_{(2)}c_{(2)})= \sigma(a_{(1)},b_{(1)})\sigma (a_{(2)}b_{(2)},c),
\end{gather*}
for $a,b,c \in {\cal Q}_{0}$ (Swedler's notation).
We can deform ${\cal Q}_{0}$ using~$\sigma$ to obtain a~new Hopf $\ast$-algebra~${\cal Q}_{0}^{\sigma}$.
Then ${\cal Q}_{0}^{\sigma}$ and $\hat{{\cal Q}}_{0\sigma}$ again form a~non degenerate pairing.
The product, $\ast$ and~$\kappa$ are changed by the following formulas whereas the coproduct remains unchanged.
For $a,b\in{\cal Q}_{0}$,
\begin{gather}
   a._{\sigma}b:= \sigma^{-1}(a_{(1)},b_{(1)})a_{(2)}b_{(2)}\sigma (a_{(3)},b_{(3)}),\nonumber
\\
 \label{eq7}  a^{*_{\sigma}}:= \sum v^{-1}(a_{(1)})a_{(2)}^{*}v(a_{(3)}),
\qquad
  \kappa_{\sigma}(a):= W(a_{(1)})\kappa(a_{(2)})W^{-1}(a_{(3)}),
\end{gather}
where
\begin{gather}\label{eq10}
  W(a)=\sigma(a_{(1)},\kappa(a_{(2)})),
\qquad
  v(a)=W^{-1}(a_{(1)})W\big(\kappa^{-1}(a_{(2)})\big).
\end{gather}
We refer the reader to~\cite[p.~65]{twist} for a~proof of co-associativity of the new product as well as other
relations of Hopf $\ast$-algebra.

We note that the unit of the deformed Hopf $\ast$-algebra remains unchanged.
For this note that
\begin{gather*}
(a._{\sigma}1)=\sigma^{-1}(a_{(1)(1)},1)a_{(1)(2)}\sigma(a_{(2)},1).
\end{gather*}
As the cocycle is normalized,
\begin{gather*}
\sigma(a,1)=\sigma^{-1}(a,1)=\epsilon(a)
\end{gather*}
for all $a\in{\cal Q}_{0}$.
So $a._{\sigma}1=\epsilon(a_{(1)(1)})a_{(1)(2)}\epsilon(a_{(2)})=a$.
Similarly it can be shown that $1._{\sigma}a=a$ for all ${\cal Q}_{0}$.

We note that ${{\cal Q}_0}^\sigma$ is a~dense Hopf $\ast$-algebra of the compact quantum group which has the discrete
quantum group $\hat{{\cal Q}_{0}}_{\sigma}$ as its dual, so we can consider the universal compact quantum group
containing ${\cal Q}_{0}^{\sigma}$ as a~dense Hopf $\ast$-algebra.
\begin{dfn}
Call the universal CQG containing ${\cal Q}_{0}^{\sigma}$ as a~dense Hopf $\ast$-algebra to be the cocycle twist of the
CQG by a~dual unitary 2-cocycle~$\sigma$ and denote it by $({\cal Q}^{\sigma},\Delta)$.
\end{dfn}

Let us now discuss how one gets a~dual unitary 2-cocycle on a~CQG from such a~dual unitary 2-cocycle on its
quantum subgroup.
Given two CQG's ${\cal Q}_{1}$, ${\cal Q}_{2}$ and a~surjective CQG morphism $\pi:{\cal Q}_{1}\rightarrow {\cal Q}_{2}$
which identif\/ies ${\cal Q}_{2}$ as a~quantum subgroup of ${\cal Q}_{1}$, i.e.~${\cal Q}_{2}\leq{\cal Q}_{1}$, it can be
shown that~$\pi$ maps the Hopf $\ast$-algebra $({\cal Q}_{1})_{0}$ onto $({\cal Q}_{2})_{0}$.
By duality we get a~map say $\hat{\pi}$ from $({\cal Q}_{2})_{0}^{\prime}$ to $({\cal Q}_{1})_{0}^{\prime}$ and it is easy to
check that this indeed maps the dense multiplier Hopf $\ast$-algebra $\hat{({\cal Q}_{2})_{0}}\subset \hat{{\cal
Q}_{2}}$ to $\hat{({\cal Q}_{1})_{0}}$.
Indeed $\hat{\pi}$ lifts to a~non degenerate $\ast$-homomorphism from ${\cal M}(\hat{{\cal Q}_{2}})$ to ${\cal
M}(\hat{{\cal Q}_{1}})$.
So given a~dual unitary 2-cocycle~$\sigma$ on ${\cal Q}_{2}$, we get a~dual unitary 2-cocycle $\sigma^{\prime}:=
(\hat{\pi}\otimes \hat{\pi})(\sigma)\in {\cal M}(\hat{{\cal Q}_{1}}\,\hat{\otimes}\,\hat{{\cal Q}_{1}})$.
It is easy to check that $\sigma^{\prime}$ is again a~dual unitary 2-cocycle on~${\cal Q}_{1}$.
We shall often use the same notation for both $\sigma^{\prime}$ and~$\sigma$, i.e.~denote $\sigma^{\prime}$ by~$\sigma$ under
slight abuse of notation for convenience.
\begin{lmma}
\label{Hopf_sub}
${\cal Q}_{2}^{\sigma}$ is a~quantum subgroup of ${\cal Q}_{1}^{\sigma^{\prime}}$.
\end{lmma}
\begin{proof}
First we claim that $\pi:({\cal Q}_{1})_{0}^{\sigma^{\prime}}\rightarrow({\cal Q}_{2})_{0}^{\sigma}$ is a~surjective Hopf
$\ast$-algebra morphism.
Since the coproducts remain unchanged, we only need to check that~$\pi$ is again a~$\ast$-algebra homomorphism.
For that observe that for $a,b\in({\cal Q}_{1})_{0}^{\sigma^{\prime}}$,
\begin{gather*}
\pi(a._{\sigma^{\prime}}b) =  \pi\big[\sigma^{\prime}(a_{(1)},b_{(1)})a_{(2)}b_{(2)}(\sigma^{\prime})^{-1}(a_{(3)},b_{(3)})\big]
\\
\hphantom{\pi(a._{\sigma^{\prime}}b)}{}
 =  \sigma^{\prime}(a_{(1)},b_{(1)})\pi(a_{(2)},b_{(2)})(\sigma^{\prime})^{-1}(a_{(3)},b_{(3)})
\\
\hphantom{\pi(a._{\sigma^{\prime}}b)}{}
 =  \sigma(\pi(a_{(1)}),\pi(b_{(1)}))\pi(a_{(2)})\pi(b_{(2)})\sigma^{-1} (\pi(a_{(3)}),\pi(b_{(3)}))
 =  \pi(a)._{\sigma}\pi(b).
\end{gather*}
Similarly we can show that $\pi(a^{\ast_{\sigma^{\prime}}})=(\pi(a))^{\ast_{\sigma}}$.
Hence ${\cal Q}_{1}^{\sigma^{\prime}}$ contains $({\cal Q}_{2})_{0}^{\sigma}$ as a~Hopf $\ast$-algebra and hence by the
universality of ${\cal Q}_{1}^{\sigma^{\prime}}$ we conclude that there is a~surjective CQG morphism from~${\cal
Q}_{1}^{\sigma^{\prime}}$ onto~${\cal Q}_{2}^{\sigma}$.
\end{proof}

\begin{lmma}
\label{back}
For a~universal CQG ${\cal Q}$ with a~dual unitary $2$-cocycle~$\sigma$, $({\cal Q}^{\sigma})^{\sigma^{-1}}\cong{\cal Q}$.
\end{lmma}

\begin{proof}
By the universality of ${\cal Q}$, it is enough to prove the Hopf $\ast$-algebra isomorphism
$({\cal Q}_{0}^{\sigma})^{\sigma^{-1}}{\cong}$ ${\cal Q}_{0}$. Again for that it is enough to check the $\ast$-algebra
structure.
For that let $a,b\in{\cal Q}_{0}$.
\begin{gather*}
a(._{\sigma})_{\sigma^{-1}}b =  \sigma^{-1}(a_{(1)},b_{(1)})a_{(2)}._{\sigma}b_{(2)}\sigma(a_{(3)},b_{(3)})
\\
\hphantom{a(._{\sigma})_{\sigma^{-1}}b}{}
 =  \sigma^{-1}(a_{(1)(1)},b_{(1)(1)})a_{(1)(2)}._{\sigma}b_{(1)(2)}\sigma(a_{(2)}, b_{(2)})
\\
\hphantom{a(._{\sigma})_{\sigma^{-1}}b}{}
 =  \sigma^{-1}(a_{(1)(1)},b_{(1)(1)})\sigma(a_{(1)(2)(1)(1)},b_{(1)(2)(1)(1)})a_{(1)(2)(1)(2)}b_{(1)(2)(1)(2)},
\\
  \sigma^{-1}(a_{(1)(2)(2)},b_{(1)(2)(2)})\sigma(a_{(2)},b_{(2)})
 =  \sigma^{-1}(a_{(1)(1)},b_{(1)(1)})\sigma(a_{(1)(2)},b_{(1)(2)})a_{(2)(1)}b_{(2)(1)},
\\
  \sigma^{-1}(a_{(2)(2)(1)},b_{(2)(2)(1)})\sigma(a_{(2)(2)(2)},b_{(2)(2)(2)})
\\
\qquad{}  =  \epsilon(a_{(1)})\epsilon(b_{(1)})a_{(2)(1)}b_{(2)(1)}\epsilon(a_{(2)(2)})\epsilon(b_{(2)(2)})
 =  \epsilon(a_{(1)})\epsilon(b_{(1)})a_{(2)}b_{(2)}
 =  ab.
\end{gather*}
Similarly we can show that the $\ast$-structure is preserved.
\end{proof}

\subsection{Unitary representations of a~twisted compact quantum group}\label{4.2}

Let ${\cal Q}$ be a~universal compact quantum group (as in the sense of Section~\ref{2.1})
with a~dual unitary 2-cocycle~$\sigma$.
In this subsection we state (with brief sketches of proofs) a~few well-known facts about the representation and the Haar
state of ${\cal Q}^{\sigma}$.
\begin{ppsn}\label{rep}\quad
\begin{enumerate}\itemsep=0pt
\item[$(i)$] For $\pi\in \operatorname{Rep}({\cal Q})$, there is an irreducible unitary representation $\pi_{\sigma}$ of ${\cal Q}^{\sigma}$,
which is the same as~$\pi$ as a~linear map.
\item[$(ii)$] $\operatorname{Rep}({\cal Q}^{\sigma})=\{\pi_{\sigma}|\pi\in \operatorname{Rep}({\cal Q})\}$.
\end{enumerate}
\end{ppsn}
\begin{proof}
(i) This follows since the $C^{\ast}$ tensor categories of unitary representations of the quantum groups ${\cal Q}$ and
${\cal Q}^{\sigma}$ are the same~\cite{bichon}.
However we give a~brief proof using the fact that both~$\hat{{\cal Q}_{0}}$ and~$(\hat{{\cal Q}_{0}})^{\sigma}$ are
$\ast$-isomorphic to $\oplus_{\pi\in \operatorname{Rep}({\cal Q})}M_{d_{\pi}}$.
Moreover recall the notation $m^{\pi}_{pq}$ and the non-degenerate pairing discussed in the Section~\ref{2.3}.
We have
\begin{gather*}
   \langle\kappa_{\sigma}(q_{ij}^{\pi})^{\ast_{\sigma}},m_{pq}^{\pi}\rangle
 =  \overline{\langle q_{ij}^{\pi},m_{pq}^{\pi\ast}\rangle}
\ \overset{\substack{\text{since the $\ast$-structure} \\ \text{does  not change}}}{=}
\  \overline{\langle q_{ij}^{\pi},m_{qp}^{\pi}\rangle}
 =  \delta_{iq}\delta_{jp}.
\end{gather*}
But we know $\langle q_{ji}^{\pi},m_{pq}^{\pi}\rangle=\delta_{iqjp}$.
Hence by non-degeneracy of the pairing $\kappa_{\sigma}(q_{ij})=q_{ji}^{\ast_{\sigma}}$, i.e.~$\pi_{\sigma}$ is again
a~unitary representation of ${\cal Q}^{\sigma}$.
\end{proof}

Thus, given a~unitary representation~$U$ (possibly inf\/inite-dimensional) of ${\cal Q}$ on
a~Hilbert space ${\cal H}$ with the spectral decomposition ${\cal H}=\oplus_{\pi\in \operatorname{Rep}({\cal Q}), i}{\cal
H}^{\pi}_{i}$, we get a~unitary representation~$U_{\sigma}$ of~${\cal Q}^{\sigma}$ on ${\cal H}$ def\/ined uniquely as
$U_{\sigma}|_{{\cal H}^{\pi}_{i}}=\pi_{\sigma}$ for all~$i$.
It is clear that $U\rightarrow U_{\sigma}$ is a~bijective correspondence between unitary representations of~${\cal Q}$
and~${\cal Q}^{\sigma}$.
\begin{ppsn}\label{mod}\quad
\begin{enumerate}\itemsep=0pt
\item[$(i)$] The Haar state for the deformed compact quantum group is the same as that of the undeformed compact quantum group.

\item[$(ii)$] The operator $F^{\sigma}_{\pi}$ corresponding to the twisted CQG ${\cal Q}^{\sigma}$ given~by
\begin{gather*}
\delta_{ik}F^{\sigma}_{\pi}(j,l)=M_{d_{\pi}}h\big(q^{\pi}_{ij}._{\sigma}q^{\pi \ast_ {\sigma}}_{kl}\big)
\end{gather*}
is related to
$F_{\pi}$ by the following
$
F^{\sigma}_{\pi}=c_{\pi}A_{\pi}^{\ast}F_{\pi}A_{\pi}$,
where $c_{\pi}$ is some positive constant and $A_{\pi}\xi:=({\rm id}\otimes v)\pi_{\sigma}\xi$, $v$~is defined by equation~\eqref{eq10}.
\end{enumerate}
\end{ppsn}
\begin{proof}
(i) This is a~straightforward consequence of the isomorphism of the coalgebra structures of ${\cal Q}$ and ${\cal
Q}^{\sigma}$.
To be more precise, denote the Haar state of ${\cal Q}$ and ${\cal Q}^{\sigma}$ by~$h$ and $h_{\sigma}$ respectively and
observing from Proposition~\ref{rep} that the matrix coef\/f\/icients of $\pi_{\sigma}$ ($\pi\in \operatorname{Rep}({\cal Q})$) are
$\{q^{\pi}_{ij}:i,j=1,\dots ,d_{\pi}\}$, we have $h_{\sigma}(q^{\pi}_{ij})=h(q^{\pi}_{ij})=0$ if~$\pi$ is non trivial and
$1$ for trivial representation.

(ii) By equation \eqref{eq1}, we see that the modular operator $\Phi|_{L^{2}(h)^{\pi}_{i}}=F^{\pi}$ for all~$\pi$ and~$i$.
Let $\Phi^{\sigma}$ be the modular operator for the CQG ${\cal Q}^{\sigma}$ and $S^{\sigma}$ be the corresponding anti
unitary operator.
Then recalling the def\/inition of the deformed $\ast$ of ${\cal Q}^{\sigma}$ (equation~\eqref{eq7}), we see that
$\Phi^{\sigma}=SC$, where $C=AB=BA$, $A$~is the operator given by $A(a)= ({\rm id}\otimes v)\Delta(a)$ and~$B$ is given by
$B(a)=(v^{-1}\otimes {\rm id})\Delta(a)$ for $a\in{\cal Q}_{0}$.
Then the modular operator $\Phi^{\sigma}=S^{\sigma\ast}S^{\sigma}=A^{\ast}B^{\ast}\Phi BA$.
We know $\Phi^{\sigma}|_{L^{2}(h)^{\pi}_{i}}=F^{\sigma}_{\pi}$ for all~$\pi$ and~$i$.
As both~$A$ and $A^{\ast}$ map ${L^{2}(h)^{\pi}_{i}}$ into itself for all~$\pi$ and~$i$, $B^{\ast}\Phi B$ also does so.
Fix some~$i$ and let $P^{\pi}_{i}$ be the projection onto $\operatorname{Sp}\{q^{\pi}_{ij}:j=1,\dots ,d_{\pi}\}$.
It is clear from the def\/inition of~$B$ that $B(q^{\pi}_{kl})=\sum\limits_{m} b^{\pi}_{km}q^{\pi}_{ml}$ and
$B^{\ast}(q^{\pi}_{kl})=\sum\limits_{m} d^{\pi}_{km}q^{\pi}_{ml}$, for some constants $b_{km}$ and $d_{km}$'s.
Thus
\begin{gather*}
   P^{\pi}_{i}\big(B^{\ast}\Phi B(q^{\pi}_{ij})\big)
 = \left(\sum\limits_{k}b_{ik}d_{ki}\right)\left(\sum\limits_{m}F_{\pi}(j,m)q^{\pi}_{im}\right),
\end{gather*}
where $\Phi(q^{\pi}_{ij})=\sum\limits_{m}F_{\pi}(j,m)q^{\pi}_{im}$.
If we denote $ \sum\limits_{k}b_{ik}d_{ki} $ by $c_{\pi,i}$, we have
\begin{gather*}
B^{\ast}\Phi B|_{L^{2}(h)^{\pi}_{i}}=c_{\pi,i}\Phi|_{L^{2}(h)^{\pi}_{i}}.
\end{gather*}
In particular, taking $c_{\pi}=c_{\pi,1}$ and denoting the restrictions of~$A$ and $A^{\ast}$ on $\operatorname{Sp}\{q^{\pi}_{1j}:1\leq j \leq d_{\pi}\}$ by $A_{\pi}$ and $A_{\pi}^{\ast}$ respectively, we write
$F^{\sigma}_{\pi}=c_{\pi}A_{\pi}^{\ast}F_{\pi}A_{\pi}$ and as $F^{\sigma}_{\pi}$ is positive invertible, $c_{\pi}$ must
be a~positive constant.
\end{proof}

\subsection{Deformation of a~von Neumann algebra by dual unitary 2-cocycles}\label{4.3}

Let ${\cal Q}$ be a~CQG with a~dual unitary 2-cocycle~$\sigma$.
Also assume that it has a~unitary representation~$V$ on a~Hilbert space ${\cal H}$ and choose a~dense subspace ${\cal
N}\subset {\cal H}$ on which~$V$ is algebraic, i.e.~$V({\cal N})\subset{\cal N}\otimes {\cal Q}_{0}$.
Then the spectral subalgebra ${\cal M}_{0}$ is SOT dense in ${\cal M}$ and ${\rm ad}_{V}({\cal M}_{0})\subset {\cal
M}_{0}\otimes {\cal Q}_{0}$ by Proposition~\ref{SOT}.
Now using the dual unitary 2-cocycle $\sigma\in {\cal M}(\hat{{\cal Q}}\,\hat{\otimes}\,\hat{{\cal Q}})$, we can def\/ine
a~new representation of ${\cal M}_{0}$ on ${\cal N}$ by
\begin{gather*}
\rho_{\sigma}(b)(\xi):= b_{(0)}\xi_{(0)}\sigma^{-1}(b_{(1)},\xi_{(1)})  \qquad \text{for} \ \ \xi\in{\cal N},
\end{gather*}
where ${\rm ad}_{V}(b)=b_{(0)}\otimes b_{(1)}$ and $V(\xi)=\xi_{(0)}\otimes \xi_{(1)}$ with Swedler's notation.
\begin{lmma}\label{new_rep}\quad
\begin{enumerate}\itemsep=0pt
\item[$(1)$] $\rho_{\sigma}(b)$ extends to an element of ${\cal B}({\cal H})$ for all $b\in{\cal M}_{0}$.

\item[$(2)$] $\rho_{\sigma}P_{\pi}$ is SOT continuous, where $P_{\pi}$ is the spectral projection corresponding to $\pi\in
\operatorname{Rep}({\cal Q})$.
\end{enumerate}
\end{lmma}
\begin{proof}
Let ${\rm ad}_{V}(b)=\sum\limits_{i=1}^{k}b_{(0)}^{i}\otimes b_{(1)}^{i}$ and
$V(\xi)=\sum\limits_{j=1}^{l}\xi_{(0)}^{j}\otimes \xi_{(1)}^{j}$.
Def\/ine $\sigma_{i}\in{\cal M}(\hat{{\cal Q}})$ by $\sigma_{i}(q_{0}):= \sigma^{-1}(b_{(1)}^{i},q_{0})$ for all
$i=1,\dots ,k$.
By Theorem~\ref{mult_rep}, we have $\Pi_{V}(\sigma_{i})\in{\cal B}({\cal H})$ for all $i=1,\dots ,k$.
By def\/inition of $\rho_{\sigma}$ on ${\cal N}$,
\begin{gather*}
\rho_{\sigma}(b)\xi= \sum\limits_{i=1}^{k}b_{(0)}^{i}\Pi_{V}(\sigma_{i})(\xi).
\end{gather*}
Using the facts that $\Pi_{V}(\sigma_{i})$ and $b_{(0)}^{i}$'s are bounded operators we can conclude that
$\rho_{\sigma}(b)\in {\cal B}({\cal H})$ for all $b\in {\cal M}_{0}$, which proves (1).

The proof of (2) is very similar and hence omitted.
\end{proof}
\begin{dfn}
We call $(\rho_{\sigma}({\cal M}_{0}))^{\prime \prime}$ the deformation of ${\cal M}$ by~$\sigma$ and denote it by
${\cal M}^\sigma$.
\end{dfn}
\begin{rmrk}\label{latest}\quad
\begin{enumerate}\itemsep=0pt
\item If there is a~weakly dense $C^*$-subalgebra ${\cal C}$ of ${\cal M}$ such that the restriction of ${\rm ad}_{V}$ leaves~${\cal C}$ invariant and also gives a~$C^*$-action on it, then the $C^*$-algebraic deformation ${\cal C}_\sigma$ in the
sense of~\cite{nesh2} can be def\/ined.
One can prove (which we plan to discuss in a~forthcoming article) that ${\cal M}_\sigma \cong ({\cal C}_\sigma)^{\prime\prime}$.

\item
It is clear from (2) of Lemma~\ref{new_rep} that $\rho_{\sigma}({\cal B})^{\prime\prime}=\rho_{\sigma}({\cal M}_{0})^{\prime\prime}$ if
${\cal B}\subset{\cal M}_{0}$ is a~SOT dense $\ast$-subalgebra of ${\cal M}_{0}$.
\end{enumerate}
\end{rmrk}

\subsection{Deformation of a~spectral triple by dual unitary 2-cocycles}

Let $({\cal A}^{\infty},{\cal H},{\cal D})$ be a~spectral triple of compact type.
Also let~$R$ be a~positive unbounded operator on ${\cal H}$ commuting with ${\cal D}$.
Let~$\sigma$ be a~dual unitary 2-cocycle on some CQG in the category ${\bf Q^{\prime}_{R}}({\cal D})$.

Given such a~dual unitary 2-cocycle we get corresponding induced dual unitary 2-cocycle on
$\widetilde{{\rm QISO}^{+}_{R}}({\cal D})$, which will again be denoted by~$\sigma$ with a~slight abuse of notation.
Let~$U$ be a~unitary representation of $\widetilde{{\rm QISO}^{+}_{R}}({\cal D})$ on ${\cal H}$.
Then~$U$ commutes with ${\cal D}$.
From now onwards we shall denote the von Neumann algebra $({\cal A}^{\infty})^{\prime\prime}$ inside ${\cal B}({\cal H})$~by
${\cal M}$.
As in Section~\ref{Section3} (right after Remark~\ref{remark3.7}),
decompose ${\cal H}$ into the eigenspaces ${\cal H}_{k}$, $k\geq 1$ and also
decompose each ${\cal H}_{k}$ into spectral subspaces, say:
\begin{gather*}
{\cal H}_{k}=\oplus_{\pi\subset {\cal I}_{k}} \mathbb{C}^{d_{\pi}}\otimes \mathbb{C}^{m_{\pi,k}},
\end{gather*}
${\cal I}_{k}\subset \operatorname{Rep}({\cal Q})$, $d_{\pi}$, $m_{\pi k}$ are as in Section~\ref{Section3}.
As $({\rm id}\otimes \phi){\rm ad}_{U}({\cal M})\subset {\cal M}$ for all bounded linear functionals~$\phi$ on
$\widetilde{{\rm QISO}^{+}_{R}}({\cal D})$, the corresponding spectral subalgebra is SOT dense.
Also for any subalgebra ${\cal A}_{0}$ of ${\cal M}$, on which ${\rm ad}_{U}$ is algebraic, we can deform it by a~dual
unitary 2-cocycle as in Section~\ref{4.3} to get a~new subalgebra in ${\cal B}({\cal H})$ which we denote by ${\cal A}_{0}^{\sigma}$.
Thus we have
\begin{thm}\label{new_sp_triple}
There is a~SOT dense $\ast$-subalgebra ${\cal A}_{0}$ in ${\cal M}$ such that
\begin{enumerate}\itemsep=0pt
\item[$(1)$] ${\rm ad}_{U}$ is algebraic over ${\cal A}_{0}$,

\item[$(2)$] $[{\cal D},a]\in {\cal B}({\cal H})$ for all $a\in {\cal A}_{0}$,

\item[$(3)$] $({\cal A}_{0}^{\sigma})^{\prime \prime}={\cal M}^{\sigma}$,

\item[$(4)$] $({\cal A}_{0}^{\sigma},{\cal H},{\cal D})$ is again a~spectral triple.
\end{enumerate}
\end{thm}
\begin{proof}
We consider the SOT-dense unital $\ast$ subalgebra ${\cal G}=\{a\in{\cal M}: [{\cal D},a]\in{\cal B}({\cal H})\}$ of
${\cal M}$.
Let $b\in{\cal G}$.
Then for a~state~$\phi$ of $\widetilde{{\rm QISO}^{+}_{R}}({\cal D})$, we have by def\/inition $({\rm id}\otimes \phi){\rm
ad}_{U}(b)\in{\cal M}$.
It also follows by the commutativity of ${\cal D}$ and~$U$ that
\begin{gather*}
[{\cal D},(({\rm id}\otimes \phi){\rm ad}_{U}(b))] =({\rm id}\otimes \phi){\rm ad}_{U}([{\cal D},b]).
\end{gather*}
Hence $({\rm id}\otimes \phi){\rm ad}_{U}(b)\in{\cal G}$ for all bounded linear functionals~$\phi$ of
$\widetilde{{\rm QISO}^{+}_{R}}({\cal D})$, i.e.~${\cal G}$ is ${\rm ad}_{U}$ invariant SOT-dense $\ast$ subalgebra of ${\cal M}$.
Now (1),(2) follow from part (3) of Lemma~\ref{new1}, taking ${\cal A}_0$ to be the span of the subspaces $P_\pi({\cal
G})$, $\pi \in \operatorname{Rep}({\cal Q})$. We also get~(3) by~2 of Remark~\ref{latest}.
To prove~(4), observe that $\forall\, a\in {\cal A}_{0}$, $\rho_{\sigma}(a)\in {\cal B}({\cal H})$ by the proof of (1) of
the Lemma~\ref{new_rep}.
So we only need to check that $[{\cal D},\rho_{\sigma}(a)]\in {\cal B}({\cal H})$ for all $a\in {\cal A}_{0}$.
With the notations used before, consider $a\in{\cal A}_{0}$ and $\xi$ in the linear span of ${\cal H}_{k}$'s, writing
${\rm ad}_{U}(a)=\sum\limits_{i=1}^{k}a_{(0)}^{i}\otimes a^{i}_{(1)}$, we have
\begin{gather*}
  [{\cal D},\rho_{\sigma}(a)](\xi)
 =  {\cal D}\rho_{\sigma}(a)(\xi)-\rho_{\sigma}(a){\cal D}(\xi)
 =  \sum\limits_{i=1}^{k}\big({\cal D} a^{i}_{(0)}\Pi_{U}(\sigma_{i})(\xi)-a_{(0)}^{i}\Pi_{U}(\sigma_{i}){\cal D}\big)(\xi)
\\
\phantom{[{\cal D},\rho_{\sigma}(a)](\xi)}{}
\overset{\substack{\text{using  the commutativity}\\ \text{of ${\cal D}$  and $U$}}}{=}  \ \sum\limits_{i=1}^{k}[{\cal D},a_{(0)}^{i}]\Pi_{U}(\sigma_{i})(\xi).
\end{gather*}
As $[{\cal D},a_{(0)}^{i}]$ is bounded for all $i=1,\dots ,k$, we conclude that $[{\cal D},\rho_{\sigma}(a)]$ is bounded
for all $a\in {\cal A}_{0}$.
\end{proof}
\begin{rmrk}
We can re-cast the deformed spectral triple $({\cal A}_{0}^{\sigma},{\cal H},{\cal D})$ in the framework of
Neshveyev--Tuset~\cite{nesh2}.
Consider the image ${\cal K}$ of the representation~$V$ viewed as an isometric Hilbert space operator from ${\cal H}$ to
${\cal H}\otimes L^{2}({\cal Q},h)$ (where~$h$ is the Haar state of ${\cal Q}$).
It can be shown that (we plan to discuss this in a~forthcoming article) $V{\cal A}_{0}^{\sigma}V^{\ast}\subset{\cal
B}({\cal H}\otimes L^{2}({\cal Q},h))$ is $\ast$-isomorphic with ${\cal A}_{0}^{\sigma}$ and $(V{\cal
A}_{0}^{\sigma}V^{\ast},{\cal K},{\cal D}\otimes 1_{L^{2}({\cal Q},h)}|_{{\cal K}})$ is the realisation of the deformed
spectral triple in the Neshveyev--Tuset framework.
\end{rmrk}

\subsection{Quantum group of orientation preserving isometries for spectral triple\\ deformed by a~dual unitary 2-cocycle}

Fix as in the previous section a~spectral triple $({\cal A}^{\infty},{\cal H},{\cal D})$ of compact type and a~positive
unbounded operator~$R$ on the Hilbert space ${\cal H}$ commuting with the Dirac operator ${\cal D}$.
Assume that there is a~dual unitary 2-cocycle~$\sigma$ on some quantum subgroup ${\cal Q}$ of
$\widetilde{{\rm QISO}^{+}_{R}}({\cal D})$.
Let~$U$,~$V$ be the unitary representations of $\widetilde{{\rm QISO}^{+}_{R}}({\cal D})$ and ${\cal Q}$ respectively on
${\cal H}$.
Let ${\cal A}_{0}$ be any SOT dense $\ast$-subalgebra of $({\cal A}^{\infty})^{\prime\prime}={\cal M}$ satisfying the conditions
of~(1) of Theorem~\ref{new_sp_triple}.
Recall from Section~\ref{4.1}, the induced dual unitary 2-cocycle $\sigma^{\prime}$ on $\widetilde{{\rm QISO}^{+}_{R}}({\cal D})$.
Since ${\rm ad}_{U}$ is algebraic over~${\cal A}_{0}$, so is ${\rm ad}_{V}$ and it is easy to see that ${\cal
A}_{0}^{\sigma^{\prime}}={\cal A}_{0}^{\sigma}$.
Now the category ${\bf Q^{\prime}_{R^{\sigma}}}({\cal A}_{0}^{\sigma},{\cal H},{\cal D})$ for some unbounded operator
$R^{\sigma}$ does not depend on the choice of the SOT dense subalgebra ${\cal A}_{0}$.
Let us abbreviate it as ${\bf Q^{\prime}_{R^{\sigma}}}({\cal D}_{\sigma})$.
We have the following:
\begin{lmma}
%\label{subobject}
$({\cal Q}^{\sigma},V_{\sigma})$ is an object in the category ${\bf Q^{\prime}_{R^{\sigma}}}({\cal D}_{\sigma})$, where
$R^{\sigma}=\Pi_{V}(v)^{\ast} R\Pi_{V}(v)$ and~$v$ is as in equation~\eqref{eq10}.
\end{lmma}

\begin{proof}
Recall the decomposition (Section~\ref{Section3}) of the Hilbert space
\begin{gather*}
{\cal H}=\oplus_{k\geq 1,\pi\in{\cal I}_{k}\subset \operatorname{Rep}({\cal Q})}\mathbb{C}^{d_{\pi}}\otimes\mathbb{C}^{m_{\pi,k}},
\end{gather*}
where $m_{\pi,k}$ is the multiplicity of the irreducible representation~$\pi$ on ${\cal H}_{k}$ and ${\cal I}_{k}$ is
some f\/inite subset of $\operatorname{Rep}({\cal Q})$.
By Theorem~\ref{rep}, $V_{\sigma}$ is again a~unitary representation of ${\cal Q}^{\sigma}$ on the Hilbert space ${\cal
H}$.
Also note that ${\rm ad}_{V}$ is algebraic over ${\cal A}_{0}$.
Let $a\in{\cal A}_{0}$ and $\xi\in {\cal N}$, where ${\cal N}$ is the subspace of ${\cal H}$ given by span of ${\cal
H}_{k}$'s.
Then we have
\begin{gather*}
V_{\sigma}(\rho_{\sigma}(a)(\xi)) =  V_{\sigma}\big(a_{(0)}\xi_{(0)}\sigma^{-1}(a_{(1)},\xi_{(1)})\big)
 =  a_{(0)(0)}\xi_{(0)(0)}\sigma^{-1}(a_{(1)},\xi_{(1)})\otimes a_{(0)(1)}\xi_{(0)(1)}
\\
\hphantom{V_{\sigma}(\rho_{\sigma}(a)(\xi))}{}
 =  a_{(0)}\xi_{(0)}\sigma^{-1}(a_{(1)(2)},\xi_{(1)(2)})\otimes a_{(1)(1)}\xi_{(1)(1)}.
\end{gather*}
On the other hand,
\begin{gather*}
  (\rho_{\sigma}\otimes [.]){\rm ad}_{V}(V_{\sigma}(\xi))
 =  (\rho_{\sigma}(a_{(0)})\xi_{(0)})\otimes a_{(1)}._{\sigma}\xi_{(1)}\\
\qquad{}
 =  a_{(0)(0)}\xi_{(0)(0)}\sigma^{-1}(a_{(0)(1)},\xi_{(0)(1)})\otimes
\sigma(a_{(1)(1)(1)},\xi_{(1)(1)(1)})a_{(1)(1)(2)}\xi_{(1)(1)(2)},
\\
  \sigma^{-1} (a_ {(1)(2)},\xi_{(1)(2)})\\
\qquad{}
 =  a_{(0)}\xi_{(0)} \sigma^{-1}(a_{(1)(1)(1)},\xi_{(1)(1)(1)})\sigma(a_{(1)(1)(2)},\xi_{(1)(1)(2)})\otimes
a_{(1)(2)(1)}\xi_{(1)(2)(1)},
\\
  \sigma^{-1}(a_{(1)(2)(2)},\xi_{(1)(2)(2)})
 =  a_{(0)}\xi_{(0)}
\epsilon(a_{(1)(1)})\epsilon(\xi_{(1)(1)})a_{(1)(2)(1)}\xi_{(1)(2)(1)}\sigma^{-1}(a_{(1)(2)(2)},\xi_{(1)(2)(2)})
\\
\qquad{} =  a_{(0)}\xi_{(0)} \epsilon(a_{(1)(1)(1)})\epsilon(\xi_{(1)(1)(1)})\otimes
a_{(1)(1)(2)}\xi_{(1)(1)(2)}\sigma^{-1}(a_{(1)(2)},\xi_{(1)(2)})
\\
\qquad{} =  a_{(0)}\xi_{(0)}\sigma^{-1}(a_{(1)(2)},\xi_{(1)(2)})\otimes a_{(1)(1)}\xi_{(1)(1)}.
\end{gather*}
So ${\rm ad}_{V_{\sigma}}(\rho_{\sigma}(a))(V_{\sigma}\xi)=V_{\sigma}(\rho_{\sigma}(a)\xi)= (\rho_{\sigma}\otimes
[.]){\rm ad}_{V}(a)(V_{\sigma}(\xi))$.
By density of ${\cal N}$ in ${\cal H}$ we conclude that
\begin{gather*}
{\rm ad}_{V_{\sigma}}(\rho_{\sigma}(a))= (\rho_{\sigma}\otimes [\cdot]){\rm ad}_{V}(a),
\end{gather*}
for all $a\in {\cal A}_{0}$.
Also for $\phi\in ({\cal Q}^{\sigma})^{\ast}$,
\begin{gather*}
({\rm id}\otimes \phi){\rm ad}_{V_{\sigma}}(\rho_{\sigma}(a))= \rho_{\sigma}(a_{(0)})\phi(a_{(1)})\in {\cal A}_{0}^{\sigma}.
\end{gather*}
So in particular $({\rm id}\otimes \phi){\rm ad}_{V_{\sigma}}({\cal A}_{0})\subset ({\cal A}_{0})^{\prime\prime}$.

$V_{\sigma}$ commutes with ${\cal D}$ since as a~linear map $V_{\sigma}$ is same as~$V$.
By Theorem~\ref{form_R},~$R$ has the form
\begin{gather*}
\oplus_{\pi\in{\cal I}_{k},k\geq 1}F_{\pi}\otimes T_{\pi,k},
\end{gather*}
for some~$T_{\pi,k}\in {\cal B}(\mathbb{C}^{m_{\pi,k}})$ so that $\Pi_{V}(v)^{\ast}R\Pi_{V}(v)$ is of the form
$\oplus_{\pi,k}A_{\pi}^{\ast}F_{\pi}A_{\pi}\otimes T_{\pi,k}$, for some~$T_{\pi,k}$ where $A_{\pi}$, $F_{\pi}$'s are as in~(ii) of Proposition~\ref{mod}.
Now by~(ii) of Proposition~\ref{mod}, we see that $F^{\sigma}_{\pi}=c_{\pi}A_{\pi}^{\ast}F_{\pi}A_{\pi}$ for some
positive constant~$c_{\pi}$.
Then $R^{\sigma}$ is of the form $ \oplus_{\pi\in{\cal I}_{k},k\geq 1}F_{\pi}^{\sigma}\otimes c_{\pi}^{-1}T_{\pi,k}$,
which implies by the `if part' of Theorem~\ref{form_R} that ${\rm ad}_{V_{\sigma}}$ preserves the $R^{\sigma}$-twisted
volume.
\end{proof}

\begin{rmrk}
%\label{withoutR}
By looking at the proof we can easily conclude that if $({\cal Q},V)$ is an object in the category ${\bf{Q}}^{\prime}({\cal
D})$, then $({\cal Q}^{\sigma},V_{\sigma})$ is an object in the category ${\bf{Q}^{\prime}({\cal D}_{\sigma})}$.
\end{rmrk}

Now replacing ${\cal Q}$ by $\widetilde{{\rm QISO}^{+}_{R}}({\cal D})$ and using the dual unitary 2-cocycle on
$\widetilde{{\rm QISO}^{+}_{R}}({\cal D})$ induced from~$\sigma$ on its quantum subgroup, we get:

\begin{crlre}
\label{qiso_sub1}
$\widetilde{{\rm QISO}^{+}_{R}}({\cal D})^{\sigma}\leq \widetilde{{\rm QISO}^{+}_{R^{\sigma}}({\cal D}_{\sigma})}$.
\end{crlre}

Thus we have the dual unitary 2-cocycle $\sigma^{-1}$ on ${\cal Q}^{\sigma}\leq
\widetilde{{\rm QISO}^{+}_{R^{\sigma}}}({\cal D}_{\sigma})$ and can deform ${\cal D}_{\sigma}$ by it.
We have the following:

\begin{lmma}
\label{pp}
$({\cal D}_{\sigma})_{\sigma^{-1}}={\cal D}$, $(R^{\sigma})^{\sigma^{-1}}=R$ and
$\widetilde{{\rm QISO}^{+}_{R^{\sigma}}}({\cal D}_{\sigma})^ {\sigma^ {-1}} \leq\widetilde{{\rm QISO}^{+}_{R}} ({\cal D})$.
\end{lmma}
\begin{proof}
Consider ${\cal M}_{1}=\operatorname{Sp} \{P_{\pi_{\sigma}}(({\cal M}_{0})^{\prime\prime}):\pi\in \operatorname{Rep}({\cal Q})\}$, where
$P_{\pi_{\sigma}}=({\rm id}\otimes \rho^{\pi_{\sigma}}){\rm ad}_{V_{\sigma}}$ (recall~$\rho^{\pi_{\sigma}}$ from Section~\ref{2.1}).
Clearly, $\rho_{\sigma^{-1}}$ is def\/ined on ${\cal M}_{1}$ and ${\cal M}_{1}$ the maximal subspace on which ${\rm
ad}_{V_{\sigma}}$ is algebraic.
As ${\rm ad}_{U_{\sigma}}$ is again a~von Neumann algebraic action of ${\cal Q}^{\sigma}$ on $({\cal
A}_{0})^{\prime\prime}=\rho_{\sigma}({\cal M}_{0})^{\prime\prime}$, we have a~SOT dense subalgebra ${\cal C}_{0}$ of $\rho_{\sigma}({\cal
M}_{0})^{\prime\prime}$ over which ${\rm ad}_{U_{\sigma}}$ is algebraic.
Then ${\rm ad}_{V_{\sigma}}$ is also algebraic over ${\cal C}_{0}$.
Hence by maximality ${\cal C}_{0}\subset{\cal M}_{1}$.
Again by SOT continuity of $\rho_{\sigma^{-1}}$ on the image of $P_{\pi_{\sigma}}$, we have $\rho_{\sigma^{-1}}({\cal
C}_{0})^{\prime\prime}=\rho_{\sigma^{-1}}({\cal M}_{1})^{\prime\prime}$.
On the other hand as ${\rm ad}_{V_{\sigma}}$ is algebraic over~$\rho_{\sigma}({\cal M}_{0})$, we have
\begin{gather*}
\rho_{\sigma}({\cal M}_{0})\subset {\cal M}_{1}
\ \Rightarrow  \ \rho_{\sigma^{-1}} (\rho_{\sigma}({\cal M}_{0}))\subset \rho_{\sigma^{-1}}({\cal M}_{1})
\
 \Rightarrow \  {\cal M}_{0}\subset\rho_{\sigma^{-1}}({\cal M}_{1}).
\end{gather*}
By maximality of ${\cal M}_{0}$, we conclude that
\begin{gather*}
\rho_{\sigma^{-1}}({\cal C}_{0})^{\prime\prime}={\cal M}_{0}^{\prime\prime}=({\cal A}^{\infty})^{\prime\prime},
\end{gather*}
which implies that ${\bf Q^{\prime}_{(R^{\sigma})^{\sigma^{-1}}}}(({\cal D}_{\sigma})_{\sigma^{-1}})={\bf Q^{\prime}_{R}} ({\cal D})$.
Also observe that $\widetilde{{\rm QISO}^{+}_{R^{\sigma}}}({\cal D}_{\sigma})^{\sigma^{-1}}$ preserves volume $\tau_{R}$ and
\begin{gather*}
{\rm ad}_{U}(a)={\rm ad}_{U_{\sigma}}(\rho_{\sigma^{-1}}(a)),
\end{gather*}
for all $a\in {\cal C}_{0}$.

But
by def\/inition $({\rm id}\otimes \phi){\rm ad}_{U_{\sigma}}(\rho_{\sigma^{-1}}(a))\subset \rho_{\sigma^{-1}}({\cal C}_{0})^{\prime\prime}=({\cal A}^{\infty})^{\prime\prime}$.
\end{proof}

Combining the above results we are in a~position to state and prove the main result of this paper.
\begin{thm}
\label{twist_main}
$\widetilde{{\rm QISO}^{+}_{R^{\sigma}}}({\cal D}_{\sigma})\cong (\widetilde{{\rm QISO}^{+}_{R}}({\cal D}))^{\sigma}$ and hence
${\rm QISO}^{+}_{R^{\sigma}}({\cal D}_{\sigma})\cong ({\rm QISO}^{+}_{R}({\cal D}))^{\sigma}$.
\end{thm}
\begin{proof}
From Lemma~\ref{pp}, we have $\widetilde{{\rm QISO}^{+}_{R^{\sigma}}}({\cal
D}_{\sigma})^{\sigma^{-1}}\leq\widetilde{{\rm QISO}^{+}_{R}} ({\cal D})$.
Then by Lemma~\ref{Hopf_sub}, we have $(\widetilde{{\rm QISO}^{+}_{R^{\sigma}}}({\cal D}_{\sigma})^{\sigma^{-1}})^{\sigma}\leq
(\widetilde{{\rm QISO}^{+}_{R}}({\cal D}))^{\sigma}$.
Since quantum isometry groups are universal, by Lemma~\ref{back}, we can conclude that
$\widetilde{{\rm QISO}^{+}_{R^{\sigma}}}({\cal D}_{\sigma})\leq (\widetilde{{\rm QISO}^{+}_{R}}({\cal D}))^{\sigma}$.
Combining this with the Corollary~\ref{qiso_sub1}, we complete the proof of the main result.
\end{proof}

Recall the category ${\bf Q^{\prime}}({\cal D})$ from Section~\ref{Section3}.
We know that in general we can not say anything about the existence of the universal object in this category.
However by looking at the proofs in this section, with the notations used in this section we have the following
\begin{crlre}
%\label{No_R}
If ${\rm QISO}^{+} ({\cal A}^{\infty},{\cal H},{\cal D})$ and ${\rm QISO}^{+} ({\cal A}_{0}^{\sigma},{\cal H},{\cal D})$ both exist,
then
\begin{gather*}
{\rm QISO}^{+} ({\cal A}^{\infty},{\cal H},{\cal D})^{\sigma}\cong {\rm QISO}^{+} ({\cal A}_{0}^{\sigma},{\cal H},{\cal D}).
\end{gather*}
\end{crlre}
\begin{rmrk}
We can also obtain similar result for the Laplacian based approach, i.e.~in the category ${\bf Q}^{\prime}_{{\cal L}}$,
whenever it makes sense.
Combining with Proposition~\ref{prevmain}, we conclude that the quantum isometry groups (${\rm QISO}^{+}_{I}$ and ${\rm QISO}^{{\cal
L}}$) for the noncommutative manifolds obtained by Rief\/fel deformation of an arbitrary compact, connected Riemannian
manifold are similar deformation of the classical isometry groups.
\end{rmrk}

\begin{rmrk}
Viewing the Rief\/fel deformation as a~special case of cocycle twist, the above theorem in fact improves the result~\cite[Theorem~5.13]{qorient} obtained by Bhowmick and Goswami by removing the assumption about a~nice dense
subalgebra on which the adjoint action of the quantum isometry group is algebraic.
In fact, techniques of this paper have enabled us to prove existence of such nice algebra in general.
\end{rmrk}

\begin{rmrk}
It should be mentioned that our result does not cover the examples constructed by Neshveyev and Tuset in~\cite{Dirac_comp}.
In that paper they considered the problem of deforming a~spectral triple coming from classical compact Lie groups using
a~Drinfeld twist.
It should be noted that although 2-cocycle is an example of a~Drinfeld twist, our scheme fails for a~general Drinfeld twist.
The main problem lies in transporting a~unitary Drinfeld twist to a~compact quantum group from its quantum subgroup.
We are grateful to Sergey Neshveyev for pointing out this dif\/f\/iculty which was overlooked mistakenly in the previous
version of this paper.
However,it seems that techniques developed in~\cite{drinfeld} can be adapted to get the main result of this paper for
the spectral triples of~\cite{Dirac_comp}.
We are working in that direction and hope to report it in another article.
\end{rmrk}

\subsection*{Acknowledgements}
Research of the f\/irst author is partially supported by Swarnajayanti Fellowship from D.S.T.\
(Govt.\
of India).
He also acknowledges kind hospitality of the Fields Institute, Toronto (Canada) where he participated in a~workshop and
conference in honour of Marc Rief\/fel.
The research of the second author is supported by CSIR, Govt.\ of India.
We are grateful to Sergey Neshveyev for pointing out a~crucial mistake in the earlier version of this paper and would
like to thank him as well as Jyotishman Bhowmick for their overall comments in general.

\pdfbookmark[1]{References}{ref}
\LastPageEnding

\end{document}